\documentclass[12pt]{amsart}
\usepackage{amsfonts}
\usepackage{amssymb,amsthm,amscd,amstext,amsmath,mathrsfs,CJK,latexsym}
\usepackage{tikz}
\usetikzlibrary{intersections,decorations.pathmorphing,decorations.fractals}

\newtheorem{theorem}{Theorem}
\newtheorem{lemma}{Lemma}

\newtheorem{cor}{Corollary}
\newtheorem{prop}{Proposition}

\newcommand{\Q}{\mathbb{Q}}
\newcommand{\C}{\mathbb{C}}

\newcommand{\CC}{\mathcal{C}}

\begin{document}

\title[]
{$\bf{|Z_{Kup}|=|Z_{Henn}|^2}$ for Lens spaces}

\author{Liang Chang}
\email{liangchang@math.ucsb.edu}
\address{Department of Mathematics\\
    University of California \\
    Santa Barbara, CA 93106\\
    U.S.A.}

\author{Zhenghan Wang}
\email{zhenghwa@microsoft.com}
\address{Microsoft Station Q\\CNSI Bldg Rm 2237\\
    University of California\\
    Santa Barbara, CA 93106-6105\\
    U.S.A.}

\thanks{}

\begin{abstract}

M.~ Hennings and G.~ Kuperberg defined quantum invariants $Z_{Henn}$ and $Z_{Kup}$ of closed oriented $3$-manifolds based on certain Hopf algebras,
respectively.  We prove that $|Z_{Kup}|=|Z_{Henn}|^2$ for lens spaces when both invariants are based on factorizable
finite dimensional ribbon Hopf algebras.

\end{abstract}

\maketitle

\section{Introduction}

A Turaev-Viro-type topological quantum field theory (TQFT) based on a spherical fusion category $\CC$ is
equivalent to the Reshetikhin-Turaev-TQFT based on the
Drinfeld center $Z(\CC)$ of $\CC$ \cite{TV2}\cite{BK}.  Consequently, $Z_{TV}(M)=|Z_{RT}(M)|^2$ for any closed oriented $3$-manifold $M$.  It is known that Hennings invariants are non-semisimple generalizations of Reshetikhin-Turaev-invariants \cite{Kerler1}, and
Kuperberg invariants are non-semisimple generalizations of Turaev-Viro-invariants \cite{BW} (in this paper, by Kuperberg invariant, we mean the one from noninvolutory Hopf algebras in \cite{Ku}).  Therefore, a similar relation might exist between the Kuperberg and Hennings invariants as first suggested in \cite{Kerler1}.

\vspace{.1in}
{\bf Problem:}
Establish a generalization of the relation between Turaev-Viro and Reshetikhin-Turaev invariants to Kuperberg and Hennings invariants.
\vspace{.1in}

One issue with the above problem is that the Kuperberg invariant $Z_{Kup}$ depends on a combing or framing of the $3$-manifold $M$, while there is
no such explicit dependence of combings or framings for the Hennings invariant $Z_{Henn}$.
Ideally, the conjectured relation would follow from a similar relation between two kinds of nonsemisimple $(2+1)$-TQFTs.  As a first step,
we prove the following relation between Kuperberg $Z_{Kup}$ and Hennings $Z_{Henn}$ invariants for the lens spaces $L(p,q)$ with $p,q\in\mathbb{N}, (p,q)=1$.

\vspace{.1in}
{\bf Main Theorem.}
{\it Let $H$ be a factorizable finite dimensional ribbon Hopf algebra and $L(p,q)$ be an oriented lens space.  Then $Z_{Kup}(L(p,q),f,H)=|Z_{Henn}(L(p,q),H)|^2$ for some suitably chosen framing $f$ of $L(p,q)$.}
\vspace{.1in}

A different choice of framing changes the Kuperberg invariant via a multiplication by a root of unity \cite{Ku}.  It follows that

\begin{cor}

$|Z_{Kup}|=|Z_{Henn}|^2$ for lens spaces and factorizable finite dimensional ribbon Hopf algebras.

\end{cor}

The contents of the paper are as follows.  In Section $2$, we recall the definitions of the Hennings and Kuperberg invariants and
set up our notations.  Finally in Section $3$, we prove the main theorem.

\section{Hennings and Kuperberg invariants}

\subsection{Some facts about Hopf algebras}

In this section, we recall some notations and structures on finite dimensional Hopf algebras.
Detail can be found in \cite{Ra1}, \cite{Ra2} and \cite{KR2}.

Let $H=(m, \Delta, S, 1, \epsilon)$ be a finite dimensional Hopf algebra over $\C$ with multiplication $m$, comultiplication $\Delta$, antipode $S$, unit $1$, and counit $\epsilon$.  We also use $1$ to denote the identity map $id$ on a Hopf algebra sometimes.

Recall that a Hopf algebra $H$ is {\it quasitriangular} if there
exists an $R$-matrix $R\in H\otimes H$.  Let $R_{ij}\in H\otimes H\otimes H$ be obtained
from $R=\sum_k s_k\otimes t_k$ by inserting the unit $1$ into the
tensor factor labeled by the index in $\{1,2,3\}\backslash \{i,j\}$.
  In a quasitriangular Hopf algebra with $R$-matrix $R=\sum_k
s_k\otimes t_k$, the special element $u=\sum_k S(t_k)s_k$
satisfies $S^2(x)=uxu^{-1}$ for $x\in H$. We use $R^{\tau}$ to denote $\sum_k t_k\otimes s_k$.

A quasitrangular Hopf
algebra is {\it ribbon} if there exists a central element
$\theta$ such that
$$\Delta(\theta)=(R^{\tau} R)^{-1}(\theta\otimes\theta),\ \ \epsilon(\theta)=1,\ \ \textrm{and} \ \ S(\theta)=\theta.$$
It can be shown that the element $G=u\theta^{-1}$ is grouplike and
$S^2(x)=GxG^{-1}$ for $x\in H$.

\subsubsection{Integrals, cointegrals and unimodular Hopf algebras}

A left integral $\lambda^L$ (respectively, right integral
$\lambda^R$) for $H$ is an element in $H^*$ which satisfies
$(id\otimes\lambda^L)\Delta(h)=\lambda^L(h)\cdot1$ (respectively,
$(\lambda^R\otimes id)\Delta(h)=\lambda^R(h)\cdot1$) for all $h\in
H$.  Dually, a left cointegral $\Lambda^L$ (respectively, right
cointegral $\Lambda^R$) for $H$ is an element in $H$ which satisfies
$h\Lambda^L=\varepsilon(h)\Lambda^L$ (respectively,
$\Lambda^Rh=\varepsilon(h)\Lambda^R$) for all $h\in H$. A Hopf
algebra $H$ is called {\it unimodular} if the space of left cointegrals
for $H$ is the same as the space of right cointegrals for $H$.

For finite dimensional Hopf algebras,
the left and right integrals (respectively, left and right
cointegrals) are unique up to scalar multiplication, and we may
choose a normalization that
$\lambda^R(\Lambda^L)=\lambda^R(S(\Lambda^L))=1$. From this, there
is an algebra homomorphism $\alpha\in H^*$, called {\it modulus} of $H$,
independent of the choice of $\Lambda^L$, such that
$\Lambda^Lh=\alpha(h)\Lambda^L$ for all $h\in H$. Likewise, there is
a grouplike element $g\in H$, called {\it comodulus} of $H$,
independent of the choice of $\lambda^R$, such that
$(id\otimes\lambda^R)\Delta(h)=\lambda^R(h)g$ for all $h\in H$. The elements $\alpha$ and $g$
are of finite order, and $\omega=\alpha(g)$ is a root of unity.

\subsubsection{Drinfeld map and factorizable Hopf algebras}

Given $Q=\sum_iQ^{(1)}_i\otimes Q^{(2)}_i\in H\otimes H$, we
define a map $f_Q:H^*\rightarrow H$ by
$f_Q(p)=\sum_ip(Q^{(1)}_i)Q^{(2)}_i$ for $p\in H^*$.  The
conditions for $R$-matrix imply that $f_R$ is an algebra
homomorphism and $f_{R^{\tau}}$ is an algebra anti-homomorphism. The map
$f_{R^{\tau}R}:H^*\rightarrow H$ is called the Drinfeld map. If the Drinfeld map for
a quasitriangular Hopf algebra $H$
is an isomorphism as a linear map of vector spaces, then $H$ is called {\it factorizable}.

\begin{prop}
If a quasitriangular Hopf algebra $H$ is factorizable, then it is unimodular.
\end{prop}

For a proof, see Prop. $3$ on page $224$ of ~\cite{Ra2}.

For a factorizable Hopf algebra $H$,
$f_{R^{\tau}R}(\lambda^R)=\Lambda^L$ and
$\lambda^R(\Lambda^L)=1$ under some normalization
(see~\cite{CW2}). This relates the left cointegral $\Lambda^L$ with
the right integral $\lambda^R$. We will use such a pair of related integral and
cointegral throughout this paper.

In this paper, we work with factorizable finite
dimensional ribbon Hopf algebras. For such Hopf algebras, we
use $\Lambda$ to denote the left and right cointegrals for
$H$. The comodulus $\alpha$ is the counit $\varepsilon$. The right
integral for $H$, denoted by $\lambda$, has the
following properties for all $x$ and $y$ in $H$~\cite{Ra1}:

\begin{enumerate}

\item  $\lambda(xy)=\lambda(S^2(y)x)$.
\item  $\lambda(gx)=\lambda(S(x))$, where $g$ is the comodulus of $H$.

\end{enumerate}
In particular, we have

\begin{lemma}
$\lambda(S^{-1}(x))=\lambda(xg)=\lambda(gx)=\lambda(S(x))$.
\end{lemma}

Such a $\lambda$ leads to a trace-like functional with the help of a square root $G$ of the comodulus $g$, i.e., $G^2=g$.
That is $tr:H\rightarrow \C$ by $tr(x)=\lambda(xG)=\lambda(Gx)$ such that $tr(xy)=tr(yx)$ and $tr(S(x))=tr(x)$ for all $x$ and $y$
in $H$.

The following lemma is important for the proof of the main theorem. Let
$\Delta^{(n-1)}(x)=\sum_{(x)} x_{(1)}\otimes\ldots\otimes x_{(n)}$
be in the Sweedler notation for iterated comultiplication. In this paper, we omit the summation
symbol, i.e.,
we write $\Delta^{(n-1)}(x)=x_{(1)}\otimes\ldots\otimes x_{(n)}$.

\begin{lemma}
For $p\in H^*$ and $n\in\mathbb{N}$, we have
\begin{eqnarray*}\Delta^{(n-1)}(f_{R^{\tau}R}(p))&=&f_{R^{\tau}}(p_{(1)})f_R(p_{(2n-1)})\otimes f_{R^{\tau}}(p_{(2)})f_R(p_{(2n-2)})\otimes\\
& &\ldots\otimes f_{R^{\tau}}(p_{(n-1)})f_R(p_{(n+1)})\otimes f_{R^{\tau}R}(p_{(n)})
\end{eqnarray*}
In particular, since $f_{R^{\tau}R}(\lambda)=\Lambda$, we have
$$\Lambda_{(k)}=f_{R^{\tau}}(\lambda_{(k)})f_R(\lambda_{(2n-k)})\ \ for\ \  k=1,\ldots,n-1$$
and $\Lambda_{(n)}=f_{R^{\tau}R}(\lambda_{(n)})$.
\end{lemma}

\begin{proof}
Recall that $f_R$ is a coalgebra antihomomorphism and $f_{R^{\tau}}$ is a
coalgebra homomorphism, i.e., for $p\in H^*$,

\begin{eqnarray*}
  & &\Delta(f_R(p))=(f_R\otimes f_R)(\Delta(p))=f_R(p_{(2)})\otimes f_R(p_{(1)}),\\
  & &\Delta(f_{R^{\tau}}(p))=(f_{R^{\tau}}\otimes f_{R^{\tau}})(\Delta^{op}(p))=f_{R^{\tau}}(p_{(1)})\otimes f_{R^{\tau}}(p_{(2)}).
\end{eqnarray*}

These two properties follow from the definition of the $R$-matrix: $(\Delta\otimes id)(R)=R_{13}R_{23}$ and $(id\otimes\Delta)(R)=R_{13}R_{12}$.\\
The proof of the lemma is by induction. When $n=2$,

\begin{eqnarray*}
  \Delta(f_{R^{\tau}R}(p))&=&\Delta(f_{R^{\tau}}(p_{(1)})f_R(p_{(2)}))\\
  &=&\Delta(f_{R^{\tau}}(p_{(1)}))\Delta(f_R(p_{(2)}))\\
  &=&f_{R^{\tau}}(p_{(1)})f_R(p_{(4)})\otimes f_{R^{\tau}}(p_{(2)})f_R(p_{(3)})\\
  &=&f_{R^{\tau}}(p_{(1)})f_R(p_{(3)})\otimes f_{R^{\tau}R}(p_{(2)})
\end{eqnarray*}

Suppose the lemma is true for $n=k$, then when $n=k+1$,

\begin{eqnarray*}
  \Delta^{(k)}(f_{R^{\tau}R}(p))&=&(id\otimes\cdots\otimes id\otimes\Delta)(\Delta^{(k-1)}(f_{R^{\tau}R}(p)))\\
  &=&(id\otimes\cdots\otimes id\otimes\Delta)(f_{R^{\tau}}(p_{(1)})f_R(p_{(2k-1)})\otimes f_{R^{\tau}}(p_{(2)})f_R(p_{(2k-2)})\otimes\\
  & &\ldots\otimes f_{R^{\tau}}(p_{(k-1)})f_R(p_{(k+1)})\otimes f_{R^{\tau}R}(p_{(k)}))\\
  &=&f_{R^{\tau}}(p_{(1)})f_R(p_{(2k-1)})\otimes f_{R^{\tau}}(p_{(2)})f_R(p_{(2k-2)})\otimes\\
  & &\ldots\otimes f_{R^{\tau}}(p_{(k-1)})f_R(p_{(k+1)})\otimes \Delta(f_{R^{\tau}R}(p_{(k)}))\\
  &=&f_{R^{\tau}}(p_{(1)})f_R(p_{(2k+1)})\otimes f_{R^{\tau}}(p_{(2)})f_R(p_{(2k)})\otimes\\
  & &\ldots\otimes f_{R^{\tau}}(p_{(k-1)})f_R(p_{(k+3)})\otimes f_{R^{\tau}}(p_{(k)})f_R(p_{(k+2)})\otimes f_{R^{\tau}R}(p_{(k+1)})
\end{eqnarray*}
Hence, the lemma holds for all $n\in\mathbb{N}$.
\end{proof}

\subsubsection{Examples}

Factorizable finite dimensional ribbon Hopf algebras include the following important examples:

\begin{enumerate}

\item $U_qsl(2,\mathbb{C})$ at an odd root of unity.
Let $q$ be an $l$-th primitive root of unity with $l$ an odd integer $\geq 3$.

$U_qsl(2,\mathbb{C})$ is generated by $E$, $F$ and $K$ with the
following relations:
$$E^l=F^l=0,\ \ K^l=1$$ and the Hopf algebra structure given by
$$KE=q^2EK,\ \ KF=q^{-2}FK,\ \ [E, F]=\frac{K-K^{-1}}{q-q^{-1}},$$
$$\Delta(E)=1\otimes E+E\otimes K,\ \ \Delta(F)=K^{-1}\otimes F+F\otimes 1, \ \ \Delta(K)=K\otimes K,$$
$$\varepsilon(E)=\varepsilon(F)=0,\ \ \varepsilon(K)=1,$$ $$S(E)=-EK^{-1},\ \ S(F)=-KF, \ \ S(K)=K^{-1}.$$
It is factorizable and ribbon with the following $R$-matrix and ribbon element
  $$R=\frac{1}{l}\sum_{0\leq m,i,j\leq l-1}\frac{(q-q^{-1})^m}{[m]!}q^{m(m-1)/2+2m(i-j)-2ij}E^mK^i\otimes F^mK^j,$$
  $$\theta=\frac{1}{l}(\sum_{s=0}^{l-1}q^{s^2})(\sum_{0\leq m,j\leq l-1}\frac{(q^{-1}-q)^m}{[m]!}q^{-\frac{1}{2}m+mj+\frac{1}{2}(j+1)^2}F^mE^mK^j)$$
Its right integral, two-sided cointegral and comodulus are
  $$\lambda(F^mE^nK^j)=\delta_{m,l-1}\delta_{n,l-1}\delta_{j,1},~~~\Lambda=F^{l-1}E^{l-1}\sum_{j=0}^{l-1}K^j,~~~g=K^2$$

\item The Drinfeld double $D(H)$ of a finite dimensional Hopf algebra $H$ is factorizable. By~\cite{KR2}, it has a ribbon element if and only if
$$S^2(h)=l( ({\beta}^{-1}\otimes \textrm{id} \otimes \beta) {\Delta}^2(h))l^{-1}$$
for all $h\in H$, where $l$ and $\beta$ are grouplike elements of $H$ and $H^*$, respectively, which satisfy $l^2=g$ and $\beta^2=\alpha$.
\end{enumerate}

\subsection{Hennings invariant}

Let $(H,R,\theta)$ be a unimodular finite dimensional ribbon Hopf algebra with $\lambda(\theta)\lambda(\theta^{-1})\neq0$.

\subsubsection{Kauffman-Radford version of the Hennings invariant}

We recall the Kauffman-Radford version of the Hennings invariant \cite{KR}.  First $(H,R,\theta)$ gives rise to a regular isotopy invariant $TR(L,H)$ for framed links $L$ as follows: given any link diagram $L_D$ of $L$, decorate each crossing of $L_D$
with the tensor factors from the
$R$-matrix $R=\sum_i s_i\otimes t_i$ as below.

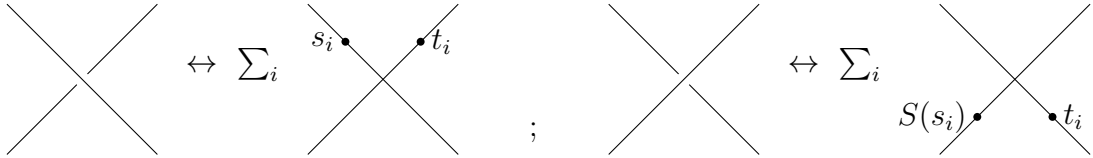
\vspace{.1in}

\begin{tikzpicture}\centering
\draw (1,1)--(-1,-1);
\draw [color=white,line width=2mm] (-1,1)--(1,-1);
\draw (-1,1)--(1,-1);
\node (=) at (2,0.2) {$\leftrightarrow~\sum_i$};
\begin{scope}[xshift = 4cm]
\draw (-1,1)--(1,-1);
\draw (1,1)--(-1,-1);
\fill[black, opacity=1] (-0.5,0.5) circle (1.5pt) node[anchor=east] {$s_i$};
\fill[black, opacity=1] (0.5,0.5) circle (1.5pt) node[anchor=west] {$t_i$};
\end{scope}
\node (=) at (6,-0.7) {$;$};
\begin{scope}[xshift = 8cm]
\draw (-1,1)--(1,-1);
\draw [color=white,line width=2mm] (-1,-1)--(1,1);
\draw (-1,-1)--(1,1);
\node (=) at (2,0.2) {$\leftrightarrow~\sum_i$};
\begin{scope}[xshift = 4.4cm]
\draw (-1,1)--(1,-1);
\draw (1,1)--(-1,-1);
\fill[black, opacity=1] (-0.5,-0.5) circle (1.5pt) node[anchor=east] {$S(s_i)$};
\fill[black, opacity=1] (0.5,-0.5) circle (1.5pt) node[anchor=west] {$t_i$};
\end{scope}
\end{scope}
\end{tikzpicture}

\vspace{.1in}

Once all the crossings of $L_D$ have been decorated, let $D_L$ be the
labeled diagram immersed in the plane, where all crossings became $4$-valent vertices. The Hopf algebra elements on $D_L$
may slide across maxima or minima of $D_L$ on the same component at the expense of the application
of the antipode or its inverse. Passing through an extremum in a clockwise direction
introduces $S^{-1}$ and passing through an extremum in a counterclockwise direction introduces $S$ as below.\\
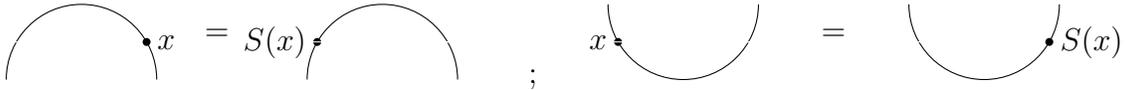
\[
\begin{tikzpicture}
  \draw [name path=circle-1] (-1,0) arc(0:180:1);
  \draw [color=white, name path=line](-3,0.5)--(3,0.5);
  \fill[black, opacity=1, name intersections={of= circle-1 and line}] (intersection-1) circle (1.5pt) node[anchor=west] {$x$};
  \node (=) at (-0.2,0.6) {=};
  \begin{scope}[xshift = 4cm]
  \draw [name path=circle-1] (-1,0) arc(0:180:1);
  \draw [color=white, name path=line](-3,0.5)--(3,0.5);
  \fill[black, opacity=1, name intersections={of= circle-1 and line}] (intersection-2) circle (1.5pt) node[anchor=east] {$S(x)$};
  \end{scope}
  \node (;) at (4,0) {$;$};
  \begin{scope}[xshift = 4cm]
  \draw [name path=circle-1] (1,1) arc(180:360:1);
  \draw [color=white, name path=line](-3,0.5)--(3,0.5);
  \fill[black, opacity=1, name intersections={of= circle-1 and line}] (intersection-1) circle (1.5pt) node[anchor=east] {$x$};
  \node (=) at (4,0.6) {=};
  \begin{scope}[xshift = 4cm]
  \draw [name path=circle-1] (1,1) arc(180:360:1);
  \draw [color=white, name path=line](-3,0.5)--(3,0.5);
  \fill[black, opacity=1, name intersections={of= circle-1 and line}] (intersection-2) circle (1.5pt) node[anchor=west] {$S(x)$};
  \end{scope}
  \end{scope}
\end{tikzpicture}
\]
To define $TR(L, H)$, slide all the Hopf algebra elements on the same
component into one vertical portion of the same component. Along a vertical
line, all the Hopf algebra elements on the same component of $D_L$ are multiplied together.\\
\[
\begin{tikzpicture}
\draw (0,1)--(0,-1);
\fill[black, opacity=1] (0,0.5) circle (1.5pt) node[anchor=west] {$y$};
\fill[black, opacity=1] (0,-0.5) circle (1.5pt) node[anchor=west] {$x$};
\node (=) at (1,0) {=};
\begin{scope}[xshift = 2cm]
\draw (0,1)--(0,-1);
\fill[black, opacity=1] (0,0) circle (1.5pt) node[anchor=west] {$xy$};
\end{scope}
\end{tikzpicture}
\]
The final juxtaposition of labeled elements at the chosen points
gives rise to a product $w_i\in H$ for the $i$-th component of $L_D$. Let
$d_i$ be the Whitney degree of this component obtained by traversing
it upward from the chosen vertical portion. The Whitney degree is
the total number of turns of the tangent vector as one traverses
the curve in the given direction. For example:\\
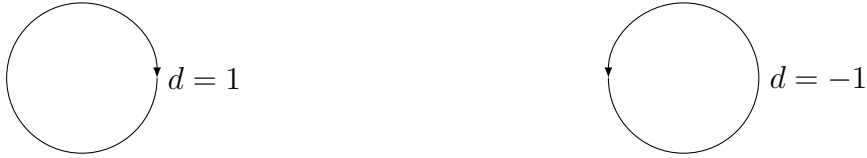
\[
\begin{tikzpicture}[>=latex]
  \draw [->](1,0) arc(360:0:1);
  \draw (1,0) node[anchor=west] {$d=1$};
  \begin{scope}[xshift=6cm]
  \draw [->](1,0) arc(-180:180:1);
  \draw (3,0) node[anchor=west] {$d=-1$};
  \end{scope}
\end{tikzpicture}
\]
Define $TR(L_D,H)=tr(w_1G^{d_1})\cdots tr(w_{c(L)}G^{d_{c(L)}})$, where $c(L)$ denotes the number of
components of $L$. This quantity is invariant under Reidemeister $II$ and $III$ moves, hence is a
regular isotopy invariant of the framed link $L$. Moreover, if
$\lambda(\theta)\lambda(\theta^{-1})\neq 0$, which is always true when $H$
is factorizable \cite{CW2}, then
\begin{equation*}
Z_{Henn}(M(L),H)=[\lambda(\theta)\lambda(\theta^{-1})]^{-\frac{c(L)}{2}}[\lambda(\theta)/\lambda(\theta^{-1})]^{-\frac{\sigma(L)}{2}}TR(L,H)\tag{2.2}
\end{equation*}
is an invariant of the closed oriented $3$-manifold $M(L)$ obtained from surgery on the framed link $L$ with
the blackboard framing, and $\sigma(L)$ denotes the signature of the
framing matrix of $L$.

\subsubsection{Properties of Hennings invariant}

Given a closed oriented manifold $M$, the symbol $\overline{M}$ denotes the same manifold with the opposite orientation.

We have the following from \cite{H}:

\begin{enumerate}

\item $Z_{Henn}(M_1\# M_2, H)=Z_{Henn}(M_1, H) Z_{Henn}(M_2, H),$

\item $Z_{Henn}(\overline{M},H)=\overline{Z_{Henn}(M,H)}$

\end{enumerate}

\subsection{Kuperberg invariant}

Let $H$ be any finite dimensional Hopf algebra.  In the following, we briefly recall some terminologies from \cite{Ku}.  For detail, see Section $2$ of \cite{Ku}.

\subsubsection{Kuperberg combings}

Given a Heegaard
diagram of a closed connected oriented $3$-manifold $M$, Kuperberg
referred to the attaching curves $c_l$'s of the $2$-handles of one handlebody as lower circles
and the attaching curves $c_u$'s of the $2$-handles of the other handlebody as upper
circles. Note that this choice is arbitrary.  A Heegaard diagram on a Heegaard surface $F$ of genus $g$ is called
{\it F-minimal} if the Heegaard diagram consists of $g$ lower circles and $g$ upper circles. In the sequel, we will simply call an $F$-minimal Heegaard diagram a minimal Heegaard diagram.  The orientation of $M$ induces an orientation on its Heegaard surface
$F$ by appending a normal vector that points from the lower side to the upper side to a positive tangent basis at a point on
$F$ which extends to a positive basis for $M$.  Define a {\it combing} on a minimal
Heegaard diagram on surface $F$ to be a vector field on $F$ with
$2g$ singularities of index $-1$, one on each circle, and one
singularity of index $+2$ disjoint from all circles. The singularity
of index $-1$ on a given circle, which is called the base point of
the circle, should not lie on a crossing and the two
outward-pointing vectors should be tangent to the circle.  Combings of Heegaard diagrams can be used to represent
combings of $3$-manifolds due to the following fact.
\begin{prop}\label{combing}
  Any combing $b$ of a minimal Heegaard diagram of $M$ can be extended to a combing
  $\overline{b}$ of $M$. Conversely, any combing of $M$ is homotopic to
  the Kuperberg extension of some combing of the minimal Heegaard
  diagram.
\end{prop}

For a proof of the proposition, see Section $2$ of \cite{Ku}.

\subsubsection{Twist front and rotation numbers}

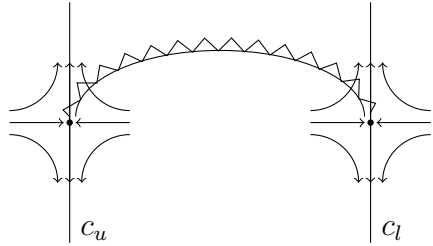
\begin{figure}\centering
  \begin{tikzpicture}[scale=0.8]
  \draw (0,-2)--(0,2);
  \fill[black, opacity=1] (0,0) circle (1.5pt);
  \draw[->] (0,0)--(0,1);
  \draw[->] (0,0)--(0,-1);
  \draw[->] (1,0)--(0.1,0);
  \draw[->] (-1,0)--(-0.1,0);
  \draw[->] (1,0.2) arc(270:180:0.8);
  \draw[->] (-1,0.2) arc(270:360:0.8);
  \draw[->] (1,-0.2) arc(90:180:0.8);
  \draw[->] (-1,-0.2) arc(90:0:0.8);
  \draw (5,-2)--(5,2);
  \fill[black, opacity=1] (5,0) circle (1.5pt);
  \draw[->] (5,0)--(5,1);
  \draw[->] (5,0)--(5,-1);
  \draw[->] (6,0)--(5.1,0);
  \draw[->] (4,0)--(4.9,0);
  \draw[->] (6,0.2) arc(270:180:0.8);
  \draw[->] (4,0.2) arc(270:360:0.8);
  \draw[->] (6,-0.2) arc(90:180:0.8);
  \draw[->] (4,-0.2) arc(90:0:0.8);
  \draw[decorate, decoration=zigzag] (0,0.1) arc(180:0:2.5 and 1.2);
  \draw (0.1,0.1) arc(180:0:2.4 and 1.1);
  \draw (0,-1.8) node[anchor=west] {$c_u$};
  \draw (5,-1.8) node[anchor=west] {$c_l$};
\end{tikzpicture}
\caption{\it Twist front}
\label{fig:twist front}
\end{figure}

Given a combing $b_1$ of $M$, by Prop. \ref{combing} we may assume it is
extended from some combing of a minimal Heegaard diagram $D$. A framing of $M$ can be obtained from another combing $b_2$ that is orthogonal to $b_1$: the third combing $b_3$ of $M$ is determined by the orientation of $M$.  To
describe such a framing $(b_1,b_2)$ of $M$, where $b_2$ is an orthogonal
combing to $b_1$, it suffices to give $b_1$ as a diagram combing
and then specify $b_2$ on the Heegaard surface $F$ and on all
upper and lower disks.  Kuperberg introduced twist fronts to encode
the position of $b_2$. A {\it twist front} is an arc along which $b_2$ is
normal to $F$ and points from the lower to the upper handlebody. A twist
front is transversely oriented in the direction that $b_2$ rotates
by the right-hand rule relative to $b_1$ and transverse orientation
is presented by the zigzag symbol as in Fig.~\ref{fig:twist front}.

To define the Kuperberg invariant, orient all Heegaard circles. Let $f=(b_1,b_2)$ be a framing
from the minimal Heegaard diagram $D$. For each point $p$ on some
circle $c$ of $D$ with base point $o$, we define $\psi(p)$ to be the
counterclockwise rotation of the tangent to $c$ relative to $b_1$
from $o$ to $p$ in units of $1=360^\circ$. If $p$ is a crossing,
then two rotation angles $\psi_l(p)$ and $\psi_u(p)$ are defined.
$\psi_l$ and $\psi_u$ are defined to be the total counterclockwise
rotation on the lower and upper circles. Let $\phi(p)$ be the total right-handed twist of $b_2$
around $b_1$ from $o$ to $p$, and similarly define $\phi_l(p)$ and
$\phi_u(p)$. Using twist fronts, we can compute $\phi(p)$ as the
total signs of all fronts crossed from $o$ to $p$, not counting the
front that terminates at $o$ itself.

\subsubsection{The Kuperberg invariant}

First we construct some special elements from the integral and
cointegral. For any integer $n$, define $\lambda_{n-\frac{1}{2}}\in
H^*$ such that $\lambda_{n-\frac{1}{2}}(x)=\lambda^R(xg^n)$ for $x\in H$ and
$\Lambda_{n-\frac{1}{2}}=(id\otimes\alpha^n)\Delta(\Lambda^R)\in H$.
Also define the tilt map $T$ to be $T(x)=(\alpha\otimes id
\otimes\alpha^{-1})\Delta^2(S^{-2}(x))$ for $x\in H$, where $g$ (or $\alpha$) are the comodulus (or modulus) of $H$.

In Kuperberg's tensor notation, the algorithm for his invariant
$Z_{Kup}(M,f,H)$ is as follows: replace each upper circle $c_u$ with
the multiplication tensor $\mu$ with one inward arrow for each crossing
and the outward arrow with $\lambda_{m}$ at the base point, with
the arrows ordered as indicated. Here $m=-\psi(c_u)$.
Replace each lower circle $c_l$ with the comultiplication tensor
$\Delta$ with an outward arrow for each crossing and the inward
arrow with $\Lambda_{n}$ at the base point, with the arrow ordered
as indicated. Here $n=\psi(c_l)$. Replace each crossing by the
tensor $\rightarrow S^aT^b\rightarrow$ where
$a=2(\psi_l(p)-\psi_u(p))-\frac{1}{2}$, $b=\phi_l(p)-\phi_u(p)$, and
$p$ is the crossing point.

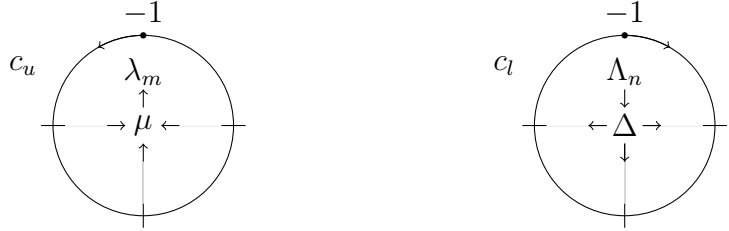
\[
\begin{tikzpicture}[scale=0.8]
  \draw (0,0) circle(1.5);
  \draw[->] (0,1.5) arc(90:120:1.5);
  \draw[->] (1.7,0)--(0.3,0);
  \draw[color=white] (1.3,0)--(0.6,0);
  \draw[->] (0,-1.7)--(0,-0.3);
  \draw[color=white] (0,-1.3)--(0,-0.6);
  \draw[->] (-1.7,0)--(-0.3,0);
  \draw[color=white] (-1.3,0)--(-0.6,0);
  \draw[->] (0,0.3)--(0,0.6);
  \draw (0,0) node {$\mu$};
  \draw (0,0.9) node {$\lambda_m$};
  \fill[black, opacity=1] (0,1.5) circle (1.5pt) node[anchor=south] {$-1$};
  \draw (-2,1) node {$c_u$};
\begin{scope}[xshift = 8cm]
  \draw (0,0) circle(1.5);
  \draw[->] (0,1.5) arc(90:60:1.5);
  \draw (1.7,0)--(0.3,0);
  \draw[color=white] (1.3,0)--(0.6,0);
  \draw[->] (0.3,0)--(0.6,0);
  \draw (0,-1.7)--(0,-0.3);
  \draw[color=white] (0,-1.3)--(0,-0.6);
  \draw[->] (0,-0.3)--(0,-0.6);
  \draw (-1.7,0)--(-0.3,0);
  \draw[color=white] (-1.3,0)--(-0.6,0);
  \draw[->] (-0.3,0)--(-0.6,0);
  \draw[->] (0,0.6)--(0,0.3);
  \draw (0,0) node {$\Delta$};
  \draw (0,0.9) node {$\Lambda_n$};
  \fill[black, opacity=1] (0,1.5) circle (1.5pt) node[anchor=south] {$-1$};
  \draw (-2,1) node {$c_l$};
\end{scope}
\end{tikzpicture}
\]

Finally, contract all tensor corresponding to circles and crossings
according to incidence.  The Kuperberg invariant
is then a big summation:
\begin{eqnarray*}
  Z_{Kup}(M,f,H)=\sum\limits_{(\Lambda)}\prod_{\substack{upper \\ circles}}\lambda\big(\cdots S^{a_i}T^{b_i}(\Lambda_{(i)})\cdots g^m\big)
\end{eqnarray*}
Here the order for multiplication and comultiplication follows
the orientations of the upper and lower circles.

For a factorizable finite dimensional ribbon Hopf algebra $H$, we have
$\alpha=\varepsilon$. So $\Lambda_{n-\frac{1}{2}}=\Lambda$ for all
integer $n$ and $T=S^{-2}$. Thus, the Kuperberg invariant is of the following
form:
\begin{equation*}\label{kuperberg invariant}
  Z_{Kup}(M,f,H)=\sum\limits_{(\Lambda)}\prod_{\substack{upper \\ circles}}\lambda\big(\cdots S^{a_i-2b_i}(\Lambda_{(i)})\cdots g^m\big)\tag{2.3}
\end{equation*}

\subsubsection{Basic properties of the Kuperberg invariant}

With suitable choices of framings, we have \cite{Ku}:

\begin{enumerate}

\item $Z_{Kup}(M_1\# M_2, H)=Z_{Kup}(M_1, H) Z_{Kup}(M_2, H),$

\item $Z_{Kup}(M, H^{*})=Z_{Kup}(\overline{M}, H^{op})=Z_{Kup}(\overline{M}, H^{cop})=Z_{Kup}(M, H).$

\end{enumerate}

\section{A relation between Kuperberg and Hennings invariants}

In this section, we prove our main theorem:

\begin{theorem}\label{Thm:Main}
  Let $H$ be a factorizable finite dimensional ribbon Hopf algebra and $L(p,q)$ be an oriented lens space, then
  $$Z_{Kup}(L(p,q),f,H)=Z_{Henn}(L(p,q)\#\overline{L(p,q)},H)$$ for
  some suitably chosen framing $f$ of $L(p,q)$.
\end{theorem}

Using $Z_{Henn}(M_1\# M_2,H)=Z_{Henn}(M_1,H)Z_{Henn}(M_2,H)$ and $Z_{Henn}(\overline{M},H)=\overline{Z_{Henn}(M,H)}$, we can deduce the version of our main theorem in the introduction.

We will calculate $Z_{Kup}(L(p,q),f,H)$ and
$Z_{Henn}(L(p,q)\#\overline{L(p,q)},H)$ through the framed Heegaard
diagram and the chain-mail link, respectively. Since $L(p,q)$ is homeomorphic to
$L(p,q+kp)$ for any integer $k$, it suffices to prove the
theorem for the case $p>q>0$.

\subsection{Chain-mail links}
Let $M$ be a closed oriented connected $3$-manifold.
We can turn a Heegaard diagram of $M$ into a surgery diagram of $M\#\overline{M}$ using the chain-mail
link introduced in~\cite{Ro}.  Let
$(F,H_1,H_2)$ be a Heegaard decomposition of $M$ with a $F$-minimal Heegaard diagram.  We will refer to $H_1$ as the lower handlebody, so
$H_2$ would be the upper handlebody.  Push the upper
circles $c_u$'s into $H_1$ slightly, then the upper circles and the
lower circles form a link in $H_1$.  All circles are framed by
thickening them into thin bands parallel to the Heegaard surface $F$. This
results in a so-called chain-mail link $C(M)\subseteq H_1$, which is in fact a
surgery presentation for $M\#\overline{M}$ (\cite{Ro}). Fig.~\ref{fig:Heegaard
L(5,2)} and Fig.~\ref{fig:chainmail L(5,2)} are the Heegaard diagram
and the corresponding chain-mail link for the Lens space $L(5,2)$, respectively.

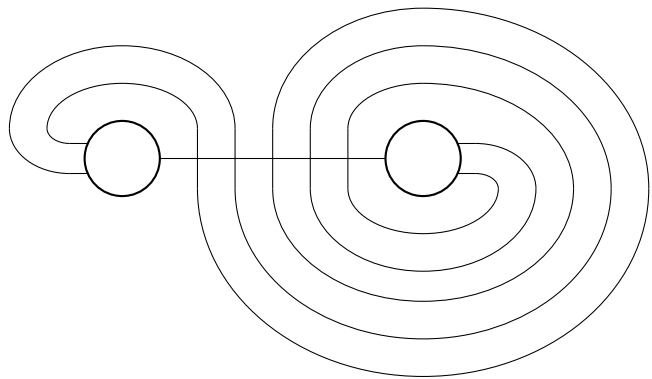
\begin{figure}\centering
\begin{tikzpicture}
  \draw [color=white, name path=hline-1] (-1.5,0.2)--(-2,0.2);
  \draw [color=white, name path=hline-2] (-1.5,-0.2)--(-2,-0.2);
  \draw [color=white, name path=hline-3] (2.5,0.2)--(3,0.2);
  \draw [color=white, name path=hline-4] (2.5,-0.2)--(3,-0.2);
  \draw [thick, name path=circle-1](-1.5,0) circle (0.5);
  \draw [thick,name path=circle-2](2.5,0) circle (0.5);
  \draw (2,0)--(-1,0);
  \draw (-0.5,-0.4)--(-0.5,0.4);
  \draw (0,-0.4)--(0,0.4);
  \draw (0.5,-0.4)--(0.5,0.4);
  \draw (1,-0.4)--(1,0.4);
  \draw (1.5,-0.4)--(1.5,0.4);
  \draw (-0.5,0.4) arc(0:180:1 and 0.6);
  \draw (-2.5,0.4) arc(180:270:0.3 and 0.2);
  \draw [name intersections={of=circle-1 and hline-1}] (-2.2,0.2)--(intersection-1);
  \draw (0,0.4) arc(0:180:1.5 and 1.1);
  \draw (-3,0.4) arc(180:270:0.8 and 0.6);
  \draw [name intersections={of=circle-1 and hline-2}] (-2.2,-0.2)--(intersection-1);
  \draw (1,-0.4) arc(180:360:1.5 and 1.1);
  \draw (4,-0.4) arc(0:90:0.8 and 0.6);
  \draw [name intersections={of=circle-2 and hline-3}] (3.2,0.2)--(intersection-1);
  \draw (1.5,-0.4) arc(180:360:1 and 0.6);
  \draw (3.5,-0.4) arc(0:90:0.3 and 0.2);
  \draw [name intersections={of=circle-2 and hline-4}] (3.2,-0.2)--(intersection-1);
  \draw (1.5,0.4) arc(180:90:1 and 0.6);
  \draw (2.5,1) arc(90:0:2 and 1.4);
  \draw (4.5,-0.4) arc(360:180:2 and 1.5);
  \draw (1,0.4) arc(180:90:1.5 and 1.1);
  \draw (2.5,1.5) arc(90:0:2.5 and 1.9);
  \draw (5,-0.4) arc(360:180:2.5 and 2);
  \draw (0.5,0.4) arc(180:90:2 and 1.6);
  \draw (2.5,2) arc(90:0:3 and 2.4);
  \draw (5.5,-0.4) arc(360:180:3 and 2.5);
\end{tikzpicture}
\caption{\it Heegaard diagram of L(5,2)}
\label{fig:Heegaard L(5,2)}
\end{figure}

\begin{figure}\centering
\begin{tikzpicture}
  \draw (-0.5,-0.4)--(-0.5,0.4);
  \draw (0,-0.4)--(0,0.4);
  \draw (0.5,-0.4)--(0.5,0.4);
  \draw (1,-0.4)--(1,0.4);
  \draw (1.5,-0.4)--(1.5,0.4);
  \draw (-0.5,0.4) arc(0:180:1 and 0.6);
  \draw (-2.5,0.4) arc(180:270:0.3 and 0.2);
  \draw (0,0.4) arc(0:180:1.5 and 1.1);
  \draw (-3,0.4) arc(180:270:0.8 and 0.6);
  \draw (1,-0.4) arc(180:360:1.5 and 1.1);
  \draw (4,-0.4) arc(0:90:0.8 and 0.6);
  \draw (1.5,-0.4) arc(180:360:1 and 0.6);
  \draw (3.5,-0.4) arc(0:90:0.3 and 0.2);
  \draw (1.5,0.4) arc(180:90:1 and 0.6);
  \draw (2.5,1) arc(90:0:2 and 1.4);
  \draw (4.5,-0.4) arc(360:180:2 and 1.5);
  \draw (1,0.4) arc(180:90:1.5 and 1.1);
  \draw (2.5,1.5) arc(90:0:2.5 and 1.9);
  \draw (5,-0.4) arc(360:180:2.5 and 2);
  \draw (0.5,0.4) arc(180:90:2 and 1.6);
  \draw (2.5,2) arc(90:0:3 and 2.4);
  \draw (5.5,-0.4) arc(360:180:3 and 2.5);
  \draw [color=white,line width=2mm] (2.3,0) arc(0:180:1.8 and 0.8);
  \draw (2.3,0) arc(0:180:1.8 and 0.8);
  \draw (-2.2,0.2) arc(270:360:0.2);
  \draw (3.2,0.2) arc(270:180:0.2);
  \draw [color=white,line width=2mm] (3,0.4) arc(0:180:2.5 and 1);
  \draw (3,0.4) arc(0:180:2.5 and 1);
  \draw (-2.2,-0.2) arc(270:360:0.2);
  \draw (3.2,-0.2) arc(270:180:0.2);
  \draw [color=white,line width=2mm] (3,0) arc(0:180:2.5 and 1.1);
  \draw (3,0) arc(0:180:2.5 and 1.1);
  \draw (2.3,0) arc(360:270:0.3);
  \draw (-1.3,0) arc(180:270:0.3);
  \draw (-1,-0.3)--(2,-0.3);
  \draw [color=white,line width=2mm](-0.5,-0.4)--(-0.5,0);
  \draw [color=white,line width=2mm](0,-0.4)--(0,0);
  \draw [color=white,line width=2mm](0.5,-0.4)--(0.5,0);
  \draw [color=white,line width=2mm](1,-0.4)--(1,0);
  \draw [color=white,line width=2mm](1.5,-0.4)--(1.5,0);
  \draw (-0.5,-0.4)--(-0.5,0);
  \draw (0,-0.4)--(0,0);
  \draw (0.5,-0.4)--(0.5,0);
  \draw (1,-0.4)--(1,0);
  \draw (1.5,-0.4)--(1.5,0);
\end{tikzpicture}
\caption{\it Chain-mail link of L(5,2)}
\label{fig:chainmail L(5,2)}
\end{figure}
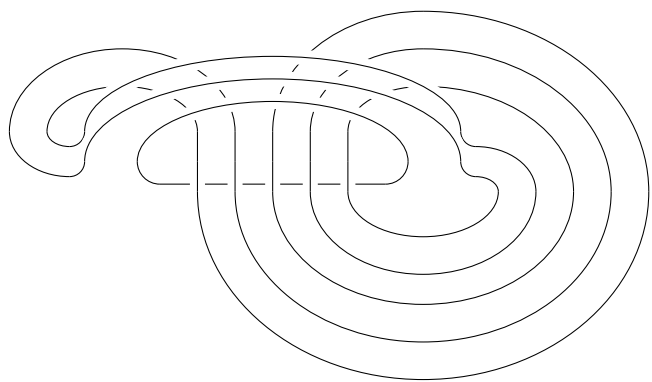

In Fig.~\ref{fig:Heegaard L(5,2)}, a $1$-handle (not drawn) is attached to the two round circles in the $2$-sphere $S^2$ regarded as
the plane together with the point at infinity.
In general, Fig.~\ref{fig:Heegaard L(p,q)} gives us a
minimal Heegaard diagram for $L(p,q)$ with $p>q>0$, where $r$ is the
remainder, i.e., $r=p-[\frac{p}{q}]q$ and $0<r<q$. Here $[x]$ means the integral part of $x$.  In the picture, the
horizontal line represents the lower circle $c_l$ (note that our lower circles are above the plane and the part of the circle over the $1$-handle is not drawn).  The upper circle
$c_u$ has $q$ strands coming out from the right circle and then
going clockwise around the right circle for $[\frac{p}{q}]-1$ times until
they meet the $q$ strands from the left circle. To make $p$
intersections with $c_l$, we let the first $r$ strands of the left $q$ strands
go around counterclockwise to match the $r$ strands of the right $q$ strands from
above. Fig.~\ref{fig:Chain-mail L(p,q)} is the corresponding
chain-mail link.

\begin{figure}\centering
\begin{tikzpicture}[scale=0.65,fill=white]
\draw [thick, name path=circle-1] (-1.5,0) circle(0.5);
\draw (1,-0.8)--(1,0.4);
\draw (1,0.4) arc(0:180:1.7 and 1.5);
\draw (-2.4,0.4) arc(180:270:0.2);
\draw [color=white, name path=hline-1] (-2,0.2)--(-1.5,0.2);
\draw [name intersections={of=circle-1 and hline-1}] (intersection-1)--(-2.2,0.2);
\draw (2,-0.8)--(2,0.4);
\draw (2,0.4) arc(0:180:2.7 and 2);
\draw (-3.4,0.4) arc(180:270:0.6);
\draw [color=white, name path=hline-3] (-2,-0.2)--(-1.5,-0.2);
\draw [name intersections={of=circle-1 and hline-3}] (intersection-1)--(-2.8,-0.2);
\draw (2,-0.8) arc(180:270:6.5 and 4.5);
\draw (8.5,-5.3) arc(270:360:6.5 and 5.3);
\draw (15,0) arc(0:90:5.5 and 3.4);
\draw (9.5,3.4) arc(90:180:5.5 and 3);
\draw (4,0.4)--(4,-0.4);
\draw [dotted,thick](4,-0.4) arc(180:220:4.5 and 3);
\draw (4,-0.4) arc(180:200:4.5 and 3);
\draw (3,-0.8) arc(180:270:5.5 and 4);
\draw (8.5,-4.8) arc(270:360:5.5 and 4.8);
\draw (14,0) arc(0:90:4.5 and 2.9);
\draw (9.5,2.9) arc(90:180:4.5 and 2.5);
\draw (5,0.4)--(5,-0.4);
\draw [dotted,thick](5,-0.4) arc(180:220:3.5 and 2.5);
\draw (5,-0.4) arc(180:200:3.5 and 2.5);
\draw [thick, name path=circle-2] (10.5,0) circle(0.5);
\draw [color=white, name path=hline-4] (11,-0.2)--(10.5,-0.2);
\draw [name intersections={of=circle-2 and hline-4}] (intersection-1)--(11.2,-0.2);
\draw (11.4,-0.4) arc(0:90:0.2);
\draw (9,-0.4) arc(180:360:1.2 and 1.2);
\draw (9,0.2)--(9,-0.4);
\draw [dotted,thick](9,0.2) arc(180:120:2 and 1);
\draw (9,0.2) arc(180:140:2 and 1);
\draw [color=white, name path=hline-6] (11,0.2)--(10.5,0.2);
\draw [name intersections={of=circle-2 and hline-6}] (intersection-1)--(11.8,0.2);
\draw (12.4,-0.4) arc(0:90:0.6);
\draw (8,-0.4) arc(180:360:2.2 and 2);
\draw (8,0.2)--(8,-0.4);
\draw [dotted,thick](8,0.2) arc(180:120:3 and 1.5);
\draw (8,0.2) arc(180:140:3 and 1.5);
\draw (3,-0.8)--(3,0.4);
\draw (3,0.4) arc(180:90:6.5 and 3.5);
\draw (9.5,3.9) arc(90:0:6.5 and 3.9);
\draw (16,0) arc(360:270:7.5 and 5.8);
\draw (8.5,-5.8) arc(270:180:7.5 and 5);
\draw (-1,0)--(5.5,0);
\draw [dashed] (5.5,0)--(7.5,0);
\draw (7.5,0)--(10,0);
\draw [dotted] (12.5,0)--(13.8,0);
\draw (1.5,0.6) node[draw,fill] { \ \ \it q \ \ };
\draw (3,0.6) node[draw,fill] {r};
\draw (2.5,-0.6) node[draw,fill] { \ \ \it q \ \ };
\draw (1,-0.6) node[draw,fill] {r};
\draw (8.5,-0.6) node[draw,fill] { \ \ \it q \ \ };
\draw (4.5,-0.6) node[draw,fill] { \ \ \it q \ \ };
\draw (-0.5,0) circle (0.1pt) node[anchor=north] {$\it c_l$};
\draw (-3.4,0.4) circle (0.1pt) node[anchor=east] {$\it c_u$};
\draw [->](1.05,-1.4)--(1.05,-1.5);
\draw [->](2.05,-1.4)--(2.05,-1.5);
\draw [->](3.05,-1.4)--(3.05,-1.5);
\draw [->](0,0)--(0.1,0);
\draw [->](10.2,-1.6)--(10.25,-1.6);
\draw [->](10.2,-2.4)--(10.25,-2.4);
\end{tikzpicture}
\caption{\it Heegaard diagram of L(p,q)} \label{fig:Heegaard L(p,q)}
\end{figure}
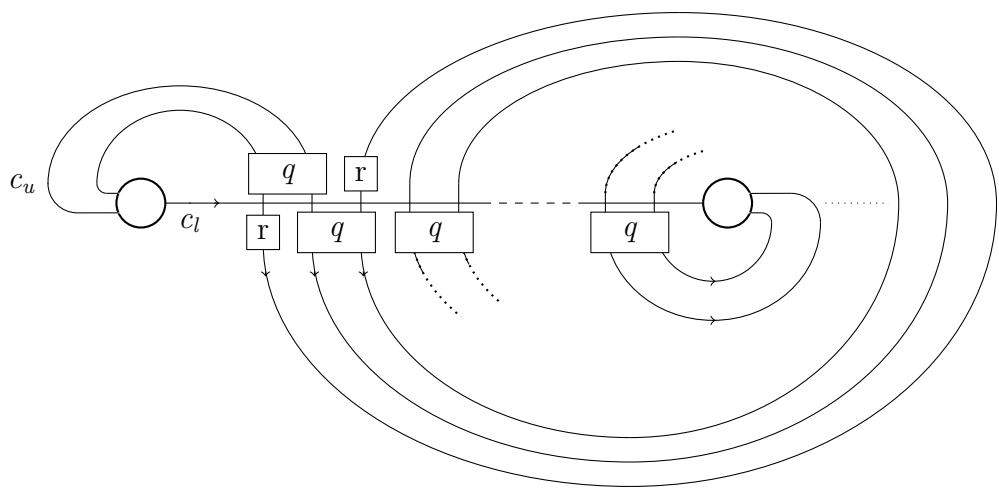

\begin{figure}\centering
\begin{tikzpicture}[scale=0.65,fill=white]
\draw (-1,-0.2) arc(180:360:5.5 and 2);
\draw (1,-0.8)--(1,0.4);
\draw (1,0.4) arc(0:180:1.7 and 1.5);
\draw (-2.4,0.4) arc(180:270:0.2);
\draw (-2,0.2)--(-2.2,0.2);
\draw (2,-0.8)--(2,0.4);
\draw (2,0.4) arc(0:180:2.7 and 2);
\draw (-3.4,0.4) arc(180:270:0.8);
\draw (-2,-0.4)--(-2.6,-0.4);
\draw [color=white,line width=2mm] (2,-0.8) arc(180:270:6.5 and 4.5);
\draw (2,-0.8) arc(180:270:6.5 and 4.5);
\draw (8.5,-5.3) arc(270:360:6.5 and 5.3);
\draw (15,0) arc(0:90:5.5 and 3.4);
\draw (9.5,3.4) arc(90:180:5.5 and 3);
\draw (4,0.4)--(4,-0.4);
\draw [color=white,line width=2mm] (4,-0.4) arc(180:230:4.5 and 3);
\draw [dotted](4,-0.4) arc(180:240:4.5 and 3);
\draw (4,-0.4) arc(180:230:4.5 and 3);
\draw [color=white,line width=2mm] (3,-0.8) arc(180:270:5.5 and 4);
\draw (3,-0.8) arc(180:270:5.5 and 4);
\draw (8.5,-4.8) arc(270:360:5.5 and 4.8);
\draw (14,0) arc(0:90:4.5 and 2.9);
\draw (9.5,2.9) arc(90:180:4.5 and 2.5);
\draw (5,0.4)--(5,-0.4);
\draw [color=white,line width=2mm] (5,-0.4) arc(180:230:3.5 and 2.5);
\draw [dotted](5,-0.4) arc(180:240:3.5 and 2.5);
\draw (5,-0.4) arc(180:230:3.5 and 2.5);
\draw (11,-0.4)--(11.2,-0.4);
\draw (11.4,-0.6) arc(0:90:0.2);
\draw [color=white,line width=2mm] (9,-0.6) arc(180:360:1.2 and 1.2);
\draw (9,-0.6) arc(180:360:1.2 and 1.2);
\draw (9,0.2)--(9,-0.6);
\draw [dotted,thick](9,0.2) arc(180:85:2 and 1);
\draw (9,0.2) arc(180:95:2 and 1);
\draw (11,0.2)--(11.8,0.2);
\draw (12.4,-0.4) arc(0:90:0.6);
\draw [color=white,line width=2mm] (8,-0.4) arc(180:360:2.2 and 2);
\draw (8,-0.4) arc(180:360:2.2 and 2);
\draw (8,0.2)--(8,-0.4);
\draw [dotted,thick](8,0.2) arc(180:85:3 and 1.5);
\draw (8,0.2) arc(180:95:3 and 1.5);
\draw (3,-0.8)--(3,0.4);
\draw (3,0.4) arc(180:90:6.5 and 3.5);
\draw (9.5,3.9) arc(90:0:6.5 and 3.9);
\draw (16,0) arc(360:270:7.5 and 5.8);
\draw [color=white,line width=2mm] (8.5,-5.8) arc(270:180:7.5 and 5);
\draw (8.5,-5.8) arc(270:180:7.5 and 5);
\draw [dotted] (5.5,0)--(7.5,0);
\draw [dotted] (12.5,0)--(13.8,0);
\draw (1.5,0.4) node[draw,fill] { \ \ \it q \ \ };
\draw (3,0.4) node[draw,fill] {r};
\draw (2.5,-0.6) node[draw,fill] { \ \ \it q \ \ };
\draw (1,-0.6) node[draw,fill] {r};
\draw (8.5,-0.6) node[draw,fill] { \ \ \it q \ \ };
\draw (4.5,-0.6) node[draw,fill] { \ \ \it q \ \ };
\draw (-2,0.2) arc(270:360:0.2);
\draw [color=white,line width=2mm] (10.8,0.4) arc(0:180:6.3 and 2.5);
\draw (10.8,0.4) arc(0:180:6.3 and 2.5);
\draw (11,0.2) arc(270:180:0.2);
\draw (-2,-0.4) arc(270:360:0.2);
\draw [color=white,line width=2mm] (10.8,-0.2) arc(0:180:6.3 and 2.5);
\draw (10.8,-0.2) arc(0:180:6.3 and 2.5);
\draw (11,0.2) arc(270:180:0.2);
\draw (11,-0.4) arc(270:180:0.2);
\draw (-2,0.2) arc(270:360:0.2);
\draw [color=white,line width=2mm] (-1,-0.2) arc(180:0:5.5 and 2);
\draw (-1,-0.2) arc(180:0:5.5 and 2);
\end{tikzpicture}
\caption{\it Chain-mail link of L(p,q)} \label{fig:Chain-mail L(p,q)}
\end{figure}
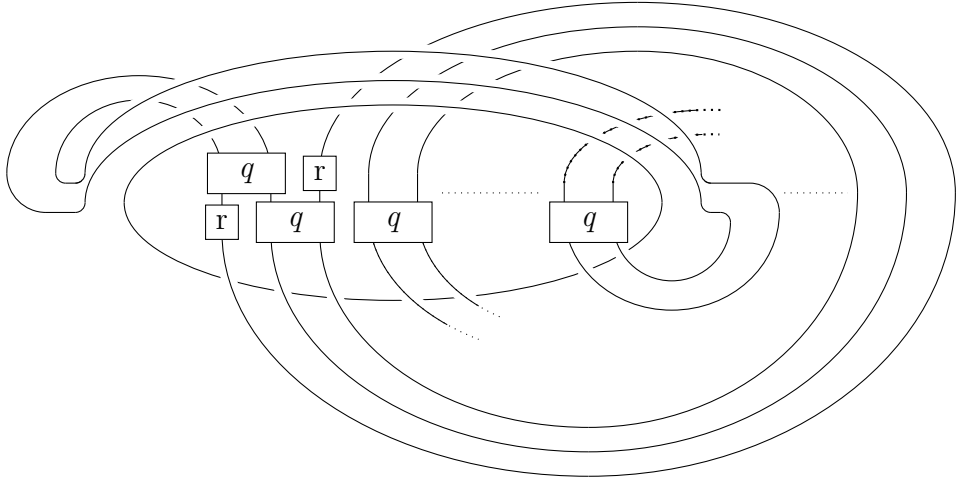

\subsection{$Z_{Kup}(L(p,q),f,H)$}
We calculate the Kuperberg invariant for $L(p,q)$
with some framing. Since the Kuperberg invariants depend on 
framings of $3$-manifolds, we need to choose a particular framing for $L(p,q)$ in order to match them with the Hennings invariants.
  The choice of framing depends on the values of $p$ and $q$.  First, let $N_1=\dfrac{q+1}{2}$ if $q$ is odd;
$N_1=\dfrac{p+q+1}{2}$ if $q$ is even. Since $p$ and $q$ are
relatively prime, $N_1$ is always a natural number. Let $N_j\equiv
N_1+(j-1)q \;\; ({\rm mod}\, p)$ such that $N_j\in \{1,\ldots,p\}$ for
$j=1,\ldots,p$.  We order the intersection points between the lower and upper circles as $1$ to $p$ from left to right in the plane, then $k_1,\cdots, k_q$ is the order of the point starting from $N_1$ to visit the first $q$ points following the orientation of the upper circle, i.e., the sequence $k_1,\cdots , k_q$ is determined by $N_{k_i}=i$ for $i=1,\cdots, q$.

In the following, we shall write $\Lambda_j=\Lambda_{(N_j)}$ as defined in Lemma $1$ for
$j=1,\ldots,p$ for short.  In other words, we rename $\Lambda_{(j)}$'s along the direction of the upper circle.  In the picture, it means the corresponding $\Lambda_{k_i}$'s are labeled on the leftmost strands.

The following technical lemma collects some symmetry properties of the indices we have just introduced.

\begin{lemma}
  (1) For $i,j\in \{1,\ldots,p\}$ such that $i+j=p+1$, we have $N_i+N_j=p+1$. As a consequence, for the two-sided cointegral $\Lambda$, we have
  $$\sum S^{\pm1}(\Lambda_{p})\otimes\cdots\otimes S^{\pm1}(\Lambda_{p+1-j})\otimes\cdots\otimes S^{\pm1}(\Lambda_{1})=\sum\Lambda_{1}\otimes\cdots\otimes \Lambda_{j}\otimes\cdots\otimes\Lambda_{p}.$$
  (2) When $q$ is odd, $k_i=1+\big[\frac{p}{q}(i-1)+\frac{q-1}{2q}\big]$ and $k_i+k_{q+2-i}=p+2$ for $i=2,\ldots,q$.\\
  (3) When $q$ is even, $k_j=1+\big[\frac{p}{q}(j-\frac{1}{2})+\frac{q-1}{2q}\big]$ and $k_j+k_{q+1-j}=p+2$ for $j=1,\ldots,q$.
\end{lemma}

\begin{proof}

Proof of $(1)$: \\
  By the definition of $N_i$, if $i+j=p+1$, then
  $$N_i+N_j\equiv 2N_1+(i+j-2)q\equiv 1 \;\; ({\rm mod}\, p)$$
  Note that $1\leq N_i,N_j\leq p$, we have $N_i+N_j=p+1$. The identity for cointegral results from the the unimodular property that $S(\Lambda)=\Lambda$. Indeed, this implies
  $$\sum S(\Lambda_{(p)})\otimes\cdots\otimes S(\Lambda_{(p+1-j)})\otimes\cdots\otimes S(\Lambda_{(1)})=\sum\Lambda_{(1)}\otimes\cdots\otimes \Lambda_{(j)}\otimes\cdots\otimes\Lambda_{(p)}$$
  Replace $\Lambda_{(N_j)}$ by $\Lambda_j$ and rearrange the order of tensor product factors, then we obtain the identity in (1).\\
  The formula for $k_i$ follows from the constraints for $N_{k_i}$:
  $$1+(i-1)p\leq N_{k_i}=N_1+(k_i-1)q\leq q+(i-1)p,$$
   if $q$ is odd, and
  $$1+ip\leq N_{k_i}=N_1+(k_i-1)q\leq q+ip,$$
   if $q$ is even.\\
   The symmetry properties of $k_i$'s are corollaries of the symmetry of the integral part function $[x]$.
   
Proof of $(2)$: \\
$$1+(i-1)p\leq N_{k_i}=N_1+(k_i-1)q\leq q+(i-1)p,$$
  then $\frac{(i-1)p}{q}\leq \frac{N_1}{q}+(k_i-1)\leq 1+\frac{(i-1)p}{q}$. So $$k_i=1+\big[\frac{p}{q}(i-1)+1-\frac{N_1}{q}\big]=1+\big[\frac{p}{q}(i-1)+\frac{q-1}{2q}\big].$$
  For $i=2,\ldots,q$, set $x_i=\frac{p}{q}(i-1)+\frac{q-1}{2q}$. Then $k_i+k_{q+2-i}=2+[x_i]+[x_{q+2-i}]$.\\
  Note that $x_i+x_{q+2-i}=p+1-\frac{1}{q}$. So
  $$[x_i]+[x_{q+2-i}]=[x_i]+\big[p+1-\frac{1}{q}-x_i\big]=p+[x_i]+\big[1-\frac{1}{q}-x_i\big].$$
  Claim: $[x_i]+\big[1-\frac{1}{q}-x_i\big]=0$.\\
  For this, let us first study the function $f(x)=[x]+[\beta-x]$ where $\beta\in[0,1)$ is a constant. It has period $T=1$ since $[x+n]=[x]$ for $n\in\mathbb{Z}$. It suffice to study it on the interval $[0,1)$. \\
  If $0\leq x\leq\beta$, $[x]+[\beta-x]=0+0=0$.\\
  If $\beta<x<1$, $[x]+[\beta-x]=0+(-1)=-1$.\\
  Hence, for $x\in\mathbb{R}$,
  \[[x]+[\beta-x]=\begin{cases}
  0&\text{for~~$x\in[n,n+\beta],~~n\in\mathbb{Z}$},\\
  -1&\text{for~~$x\in(n+\beta,n+1),~~n\in\mathbb{Z}$}.
  \end{cases}\]
  Let $\{x\}:=x-[x]$ be the fractional part of $x$. Then $[x]+[\beta-x]=0$ if and only if $\{x\}\leq\beta$.
Thus our claim is equivalent to $\{x_i\}\leq1-\frac{1}{q}$. In the following, we calculate $\{x_i\}$ case by case.\\
Let $r=(i-1)p-q\big[\frac{(i-1)p}{q}\big]$ be the remainder. Then $$\{x_i\}=\big\{\big[\frac{(i-1)p}{q}\big]+\frac{r}{q}+\frac{q-1}{2q}\big\}=\big\{\frac{r}{q}+\frac{q-1}{2q}\big\}.$$
Case 1: $1\leq r\leq\frac{q-1}{2}$.
Because $\frac{r}{q}+\frac{q-1}{2q}\leq\frac{q-1}{2q}+\frac{q-1}{2q}=1-\frac{1}{q}\leq1$, we have
$$\big\{\frac{r}{q}+\frac{q-1}{2q}\big\}=\frac{r}{q}+\frac{q-1}{2q}\leq1-\frac{1}{q}.$$
Case 2: $\frac{q+1}{2}\leq r\leq q-1$.
Because $2>\frac{r}{q}+\frac{q-1}{2q}\geq\frac{q+1}{2q}+\frac{q-1}{2q}=1$, we have
$$\big\{\frac{r}{q}+\frac{q-1}{2q}\big\}=\frac{r}{q}+\frac{q-1}{2q}-1<\frac{q}{q}+\frac{q-1}{2q}-1=\frac{1}{2}\big(1-\frac{1}{q}\big)\leq1-\frac{1}{q}.$$

Proof of $(3)$: \\
$$1+jp\leq N_{k_j}=N_1+(k_j-1)q\leq q+jp,$$
  then $\frac{jp}{q}\leq \frac{N_1}{q}+(k_j-1)\leq 1+\frac{jp}{q}$. So $$k_j=1+\big[\frac{jp}{q}+1-\frac{N_1}{q}\big]=1+\big[\frac{p}{q}(j-\frac{1}{2})+\frac{q-1}{2q}\big].$$
  Similarly as above, we set $y_j=\frac{p}{q}(j-\frac{1}{2})+\frac{q-1}{2q}$ for $j=1,\dots,q$. Then $y_j+y_{q+1-j}=p+1-\frac{1}{q}$. Thus $k_j+k_{q+1-j}=p+2$ is equivalent to
  $$[y_j]+\big[1-\frac{1}{q}-y_j\big]=0.$$
  Further more, this is equivalent to $\{y_j\}\leq1-\frac{1}{q}$.\\
  Let $r=(2i-1)p-2q\big[\frac{(2i-1)p}{2q}\big]$ be the remainder. Note that $r\neq q$ for $(p,q)=1$. Then $$\{y_j\}=\big\{\big[\frac{(2i-1)p}{2q}\big]+\frac{r}{2q}+\frac{q-1}{2q}\big\}=\big\{\frac{r}{2q}+\frac{q-1}{2q}\big\}.$$
  Case 1: $1\leq r\leq q-1$.
Because $\frac{r}{2q}+\frac{q-1}{2q}\leq\frac{q-1}{2q}+\frac{q-1}{2q}=1-\frac{1}{q}<1$, we have
$$\big\{\frac{r}{q}+\frac{q-1}{2q}\big\}=\frac{r}{2q}+\frac{q-1}{2q}\leq1-\frac{1}{q}.$$
Case 2: $q+1\leq r\leq 2q-1$.
Because $2>\frac{r}{2q}+\frac{q-1}{2q}\geq\frac{q+1}{2q}+\frac{q-1}{2q}=1$, we have
$$\big\{\frac{r}{2q}+\frac{q-1}{2q}\big\}=\frac{r}{2q}+\frac{q-1}{2q}-1<\frac{2q}{2q}+\frac{q-1}{2q}-1=\frac{1}{2}\big(1-\frac{1}{q}\big)<1-\frac{1}{q}.$$
\end{proof}

In addition, we define $k_0=1$ and $k_{q+1}=p+1$ for future use.

We set up a framed Heegaard diagram for $L(p,q)$ shown in
Fig.~\ref{fig:odd framed Heegaard} for $q$ odd, and Fig.~\ref{fig:even framed
Heegaard} for $q$ even. The framing $f$ is represented by the dashed flows and the
twist front. The twist fronts vary depending on whether $q$ is odd
or even. Two index $-1$ singularities are located at the two ends of
the twist front on the horizontal line. The right $-1$ singularity is at
the $N_1$-th intersection of the horizontal line and the upper
circle $c_u$. The lower circle $c_l$ is represented by the horizontal
line with the base point the left $-1$ singularity and oriented towards
right from the base point. To avoid the right singularity, we make $c_l$
turn around slightly when it meets this singularity. Likewise, the
upper circle $c_u$ with the base point the right singularity is oriented
towards down from its base point. The index  $+2$ singularity is
located at the infinity. The orientation of circles is shown in Fig.\ref{fig:Heegaard L(p,q)}.

\subsection{Examples}
Before we turn to the general calculation, it is helpful to examine two concrete examples.

\subsubsection{$L(2,1)$}

\begin{figure}\centering
\begin{tikzpicture}[scale=0.65,fill=white]\scriptsize
\draw [name path=circle-1] (-1.5,0) circle(0.5);
\draw (2.2,0.4)--(2.2,-0.4);
\draw (2.2,0.4) arc(0:180:2.55 and 2);
\draw (-2.9,0.4) arc(180:270:0.4);
\draw (-2,0)--(-2.5,0);
\draw (7,0)--(7.5,0);
\draw (7.5,0) arc(90:0:0.4);
\draw (7.9,-0.4) arc(360:180:1.95 and 1.5);
\draw (4,0.4)--(4,-0.4);
\draw (4,0.4) arc(180:90:2.5 and 2);
\draw (6.5,2.4) arc(90:0:2.5 and 2.4);
\draw (2.2,-0.4) arc(180:270:3.4 and 2.5);
\draw (5.6,-2.9) arc(270:360:3.4 and 2.9);
\draw [name path=circle-2] (6.5,0) circle(0.5);
\draw (-1,0)--(1.9,0);
\draw (2.5,0)--(6,0);
\draw (1.9,0) arc(180:360:0.3);
\draw [dashed](0,-4)--(0,-3);
\draw [>-<, dashed] (0,-3)--(0,3);
\draw [dashed](0,3)--(0,4);
\draw [<->, dashed] (-0.5,0)--(0.5,0);
\draw [dashed] (-1.5,0) arc(270:360:1.2 and 4);
\draw [->,dashed] (-0.3,4) arc(360:330:1.2 and 4);
\draw [dashed] (-1.5,0) arc(90:0:1.2 and 4);
\draw [->,dashed] (-0.3,-4) arc(0:40:1.2 and 4);
\draw [dashed] (-3.5,3)--(-1.5,0);
\draw [->,dashed] (-3.5,3)--(-2.5,1.5);
\draw [dashed] (-3.5,-3)--(-1.5,0);
\draw [->,dashed] (-3.5,-3)--(-2.5,-1.5);
\draw [dashed] (-3,0)--(-1.5,0);
\draw [->,dashed] (-3.5,0)--(-2.9,0);
\fill [white] (-1.5,0) circle(0.5);
\draw [thick] (-1.5,0) circle(0.5);
\draw [dashed](2.2,-4)--(2.2,-2);
\draw [<->, dashed] (2.2,-3)--(2.2,3);
\draw [dashed](2.2,2)--(2.2,4);
\draw [>-<, dashed] (1.5,0)--(2.9,0);
\draw [dashed] (0.3,4) arc(180:360:0.8 and 3);
\draw [->,dashed] (0.3,4) arc(180:210:0.8 and 3);
\draw [->,dashed] (0.3,4) arc(180:330:0.8 and 3);
\draw [dashed] (0.3,-4) arc(180:0:0.8 and 3);
\draw [->,dashed] (0.3,-4) arc(180:150:0.8 and 3);
\draw [->,dashed] (0.3,-4) arc(180:30:0.8 and 3);
\draw [dashed] (6.5,0.3) arc (270:180:4 and 3.7);
\draw [<->,dashed] (6.5,0.3) arc (270:220:4 and 3.7);
\draw [dashed] (6.5,-0.3) arc (90:180:4 and 3.7);
\draw [<->,dashed] (6.5,-0.3) arc (90:140:4 and 3.7);
\draw [dashed] (6.5,0)--(10.5,4);
\draw [<->,dashed] (6.5,0)--(8.5,2);
\draw [dashed] (6.5,0)--(10.5,-4);
\draw [<->,dashed] (6.5,0)--(8.5,-2);
\fill [white] (6.5,0) circle(0.5);
\draw [thick] (6.5,0) circle(0.5);
\draw [decorate, decoration=zigzag](0,0) arc(180:0:1.1 and 0.7);
\draw (0,0) arc(180:0:1.1 and 0.5);
\node (lambda) at (2.8,-0.6) {$\it \Lambda_{(1)}$};
\node (lambda) at (4.5,-0.4) {$\it \Lambda_{(2)}$};
\end{tikzpicture}
\caption{\it Heegaard diagram of L(2,1)}
\label{fig:Heegaard diagram L(2,1)}
\end{figure}
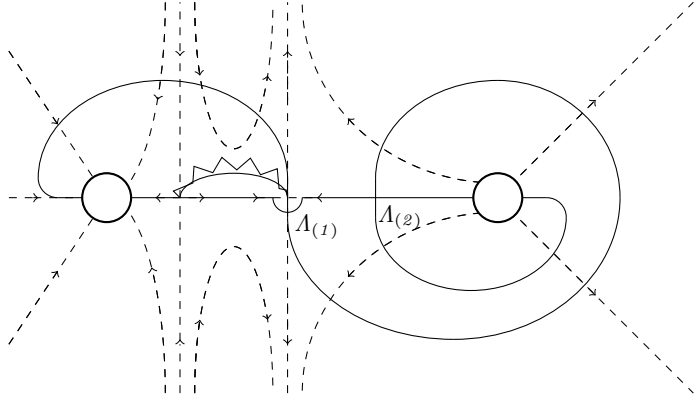
We choose a framing for $L(2,1)$ as shown in Fig.\ref{fig:Heegaard diagram L(2,1)} then go downwards from $\Lambda_{(1)}$. The power for $S$ for $\Lambda_{(1)}$ and $\Lambda_{(2)}$ are, respectively,
\begin{eqnarray*}
  &2\big(-\frac{1}{4}-0\big)-\frac{1}{2}=-1;\\
  &2\big(-\frac{1}{2}-(-\frac{1}{4})\big)-\frac{1}{2}=-1;
\end{eqnarray*}
The total rotation along the upper circle is $-\frac{1}{4}+\frac{1}{2}+\frac{1}{4}=\frac{1}{2}$ and the power of $g$ is $-\frac{1}{2}+\frac{1}{2}=0$. So the Kuperberg invariant for $L(2,1)$ is
\begin{eqnarray*}
Z_{Kup}(L(2,1))=\lambda\big(S^{-1}(\Lambda_{(1)})S^{-1}(\Lambda_{(2)})\big)=\lambda\big(\Lambda_{(2)}\Lambda_{(1)}\big)=Tr(S^{-1})
\end{eqnarray*}
where $Tr$ is the trace for vector spaces. The last equality follows from a trace formula in terms of integrals in \cite{Ra1}.

\subsubsection{$L(5,2)$}

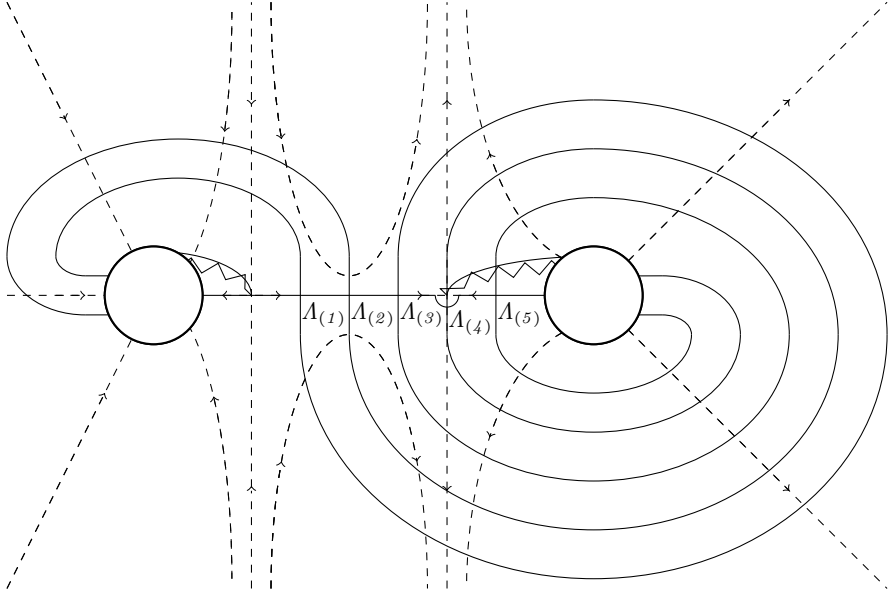
\begin{figure}\centering
\begin{tikzpicture}[scale=1.3]\scriptsize
  \draw [color=white, name path=hline-1] (-2,0.2)--(-2.5,0.2);
  \draw [color=white, name path=hline-2] (-2,-0.2)--(-2.5,-0.2);
  \draw [color=white, name path=hline-3] (2.5,0.2)--(3,0.2);
  \draw [color=white, name path=hline-4] (2.5,-0.2)--(3,-0.2);
  \draw [thick, name path=circle-1](-2,0) circle (0.5);
  \draw [thick, name path=circle-2](2.5,0) circle (0.5);
  \draw (0.88,0)--(-1.5,0);
  \draw (2,0)--(1.12,0);
  \draw (0.88,0) arc(180:360:0.12);
  \draw (-0.5,-0.4)--(-0.5,0.4);
  \draw (0,-0.4)--(0,0.4);
  \draw (0.5,-0.4)--(0.5,0.4);
  \draw (1,-0.4)--(1,0.4);
  \draw (1.5,-0.4)--(1.5,0.4);
  \draw (-0.5,0.4) arc(0:180:1.25 and 0.8);
  \draw (-3,0.4) arc(180:270:0.3 and 0.2);
  \draw [name intersections={of=circle-1 and hline-1}] (-2.7,0.2)--(intersection-1);
  \draw (0,0.4) arc(0:180:1.75 and 1.2);
  \draw (-3.5,0.4) arc(180:270:0.8 and 0.6);
  \draw [name intersections={of=circle-1 and hline-2}] (-2.7,-0.2)--(intersection-1);
  \draw (1,-0.4) arc(180:360:1.5 and 1);
  \draw (4,-0.4) arc(0:90:0.8 and 0.6);
  \draw [name intersections={of=circle-2 and hline-3}] (3.2,0.2)--(intersection-1);
  \draw (1.5,-0.4) arc(180:360:1 and 0.6);
  \draw (3.5,-0.4) arc(0:90:0.3 and 0.2);
  \draw [name intersections={of=circle-2 and hline-4}] (3.2,-0.2)--(intersection-1);
  \draw (1.5,0.4) arc(180:90:1 and 0.6);
  \draw (2.5,1) arc(90:0:2 and 1.4);
  \draw (4.5,-0.4) arc(360:180:2 and 1.5);
  \draw (1,0.4) arc(180:90:1.5 and 1.1);
  \draw (2.5,1.5) arc(90:0:2.5 and 1.9);
  \draw (5,-0.4) arc(360:180:2.5 and 2);
  \draw (0.5,0.4) arc(180:90:2 and 1.6);
  \draw (2.5,2) arc(90:0:3 and 2.4);
  \draw (5.5,-0.4) arc(360:180:3 and 2.5);
\draw [decorate, decoration=zigzag](-1,0) arc(0:90:1 and 0.34);
\draw (-1,0) arc(0:90:1 and 0.45);
\draw [decorate, decoration=zigzag](1,0) arc(180:90:1.5 and 0.3);
\draw (1,0) arc(180:90:1.5 and 0.4);
\draw [dashed](-1,-3)--(-1,-2);
\draw [>-<, dashed] (-1,-2)--(-1,2);
\draw [dashed](-1,2)--(-1,3);
\draw [<->, dashed] (-1.3,0)--(-0.7,0);
\draw [dashed] (-2,0) arc(270:360:0.8 and 2.9);
\draw [->,dashed] (-1.2,2.9) arc(360:335:0.8 and 2.9);
\draw [dashed] (-2,0) arc(90:0:0.8 and 2.9);
\draw [->,dashed] (-1.2,-2.9) arc(0:40:0.8 and 2.9);
\draw [dashed] (-3.5,3)--(-2,0);
\draw [->,dashed] (-3.5,3)--(-2.9,1.8);
\draw [dashed] (-3.5,-3)--(-2,0);
\draw [->,dashed] (-3.5,-3)--(-2.5,-1);
\draw [dashed] (-3,0)--(-1.5,0);
\draw [->,dashed] (-3.5,0)--(-2.7,0);
\draw [dashed](1,-3)--(1,-2);
\draw [<->, dashed] (1,-2)--(1,2);
\draw [dashed](1,2)--(1,3);
\draw [>-<, dashed] (0.7,0)--(1.3,0);
\draw [dashed] (-0.8,3) arc(180:360:0.8 and 2.8);
\draw [->,dashed] (-0.8,3) arc(180:210:0.8 and 2.8);
\draw [->,dashed] (-0.8,3) arc(180:330:0.8 and 2.8);
\draw [dashed] (-0.8,-3) arc(180:0:0.8 and 2.6);
\draw [->,dashed] (-0.8,-3) arc(180:150:0.8 and 2.6);
\draw [->,dashed] (-0.8,-3) arc(180:30:0.8 and 2.6);
\draw [dashed] (2.5,0.3) arc (270:180:1.3 and 2.7);
\draw [<->,dashed] (2.5,0.3) arc (270:215:1.3 and 2.7);
\draw [dashed] (2.5,-0.3) arc (90:180:1.3 and 2.7);
\draw [<->,dashed] (2.5,-0.3) arc (90:145:1.3 and 2.7);
\draw [dashed] (2.5,0)--(5.5,3);
\draw [<->,dashed] (2.5,0)--(4.5,2);
\draw [dashed] (2.5,0)--(5.5,-3);
\draw [<->,dashed] (2.5,0)--(4.5,-2);
\fill [white] (-2,0) circle(0.5);
\draw [thick] (-2,0) circle(0.5);
\fill [white] (2.5,0) circle(0.5);
\draw [thick] (2.5,0) circle(0.5);
  \node (lambda) at (-0.26,-0.2) {$\it \Lambda_{(1)}$};
  \node (lambda) at (0.24,-0.2) {$\it \Lambda_{(2)}$};
  \node (lambda) at (0.74,-0.2) {$\it \Lambda_{(3)}$};
  \node (lambda) at (1.24,-0.28) {$\it \Lambda_{(4)}$};
  \node (lambda) at (1.74,-0.2) {$\it \Lambda_{(5)}$};
\end{tikzpicture}
\caption{\it Heegaard diagram of L(5,2)}
\label{fig:Heegaard diagram L(5,2)}
\end{figure}

In Fig.\ref{fig:Heegaard diagram L(5,2)}, a framing is set up for $L(5,2)$. We start from $\Lambda_{(4)}$ and go downwards. The power for $S$ for $\Lambda_{(4)}$, $\Lambda_{(1)}$, $\Lambda_{(3)}$, $\Lambda_{(5)}$ and $\Lambda_{(2)}$ are, respectively,
\begin{eqnarray*}
  &2\big(-\frac{1}{4}-0\big)-\frac{1}{2}=-1;\\
  &2\big(0-(-\frac{1}{4}+\frac{1}{2})\big)-\frac{1}{2}=-1;\\
  &2\big(0-(-\frac{1}{4}+\frac{1}{2}-1)\big)-\frac{1}{2}=1;\\
  &2\big(-\frac{1}{2}-(-\frac{1}{4}+\frac{1}{2}-1-\frac{1}{2})\big)-\frac{1}{2}-2=3;\\
  &2\big(0-(-\frac{1}{4}+\frac{1}{2}-1-\frac{1}{2}+\frac{1}{2})\big)-\frac{1}{2}-2=3;
\end{eqnarray*}
In the last two equations, the $-2$'s result from crossing the twist front before $\Lambda_{(5)}$ and $\Lambda_{(2)}$.
The total rotation along the upper circle is $-\frac{1}{4}+\frac{1}{2}-1-\frac{1}{2}+\frac{1}{2}-\frac{3}{4}=-\frac{3}{2}$ and the power of $g$ is $\frac{3}{2}+\frac{1}{2}=2$. From this data, the Kuperberg invariant for $L(5,2)$ is
\begin{eqnarray*}
Z_{Kup}(L(5,2))&=&\lambda\big(S^{-1}(\Lambda_{(4)})S^{-1}(\Lambda_{(1)})S^{}(\Lambda_{(3)})S^{3}(\Lambda_{(5)})S^{3}(\Lambda_{(2)})g^{2}\big)\\
&=&\lambda\big(\Lambda_{(2)}\Lambda_{(5)}S^{2}(\Lambda_{(3)})S^{4}(\Lambda_{(1)})S^{4}(\Lambda_{(4)})g^{2}\big)\\
&=&\lambda\big(S^{-4}(\Lambda_{(2)})S^{-4}(\Lambda_{(5)})S^{-2}(\Lambda_{(3)})\Lambda_{(1)}\Lambda_{(4)}g^{2}\big)
\end{eqnarray*}
In this case, $N_1=4$ and $\Lambda_1$, $\Lambda_2$, $\Lambda_3$, $\Lambda_4$, $\Lambda_5$ are, respectively, $\Lambda_{(4)}$, $\Lambda_{(1)}$, $\Lambda_{(3)}$, $\Lambda_{(5)}$, $\Lambda_{(2)}$. Therefore, $k_1=2$ and $k_2=5$.

\subsection{General Calculation for $Z_{Kup}(L(p,q),f,H)$}

\subsubsection{$Z_{Kup}(L(p,q),f,H)$ when q is odd}

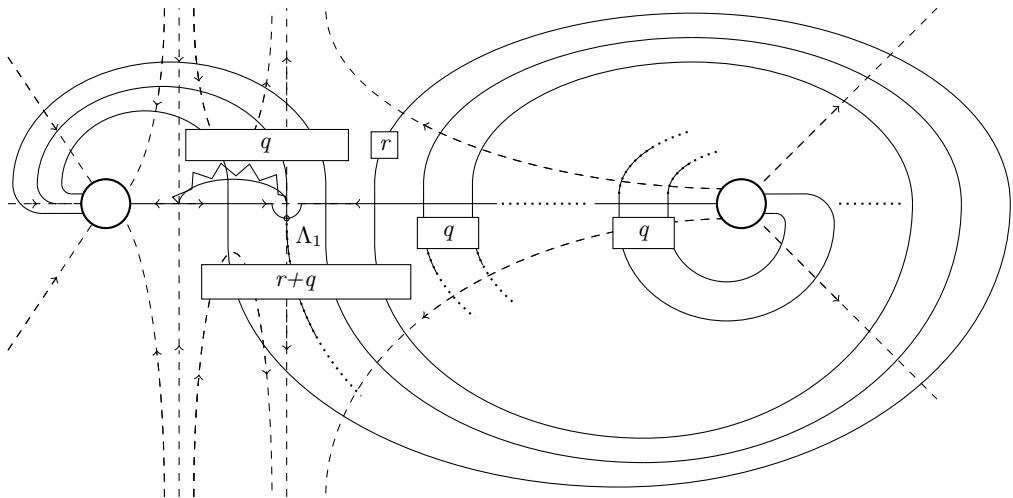
\begin{figure}\centering
\begin{tikzpicture}[scale=0.65,fill=white]\scriptsize
\draw [name path=circle-1] (-1.5,0) circle(0.5);
\draw (1,-0.8)--(1,0.4);
\draw (1,0.4) arc(0:180:1.7 and 1.5);
\draw (-2.4,0.4) arc(180:270:0.2);
\draw [color=white, name path=hline-1] (-2,0.2)--(-1.5,0.2);
\draw [name intersections={of=circle-1 and hline-1}] (intersection-1)--(-2.2,0.2);
\draw (3,-0.8)--(3,0.4);
\draw (3,0.4) arc(0:180:3.2 and 2.5);
\draw (-3.4,0.4) arc(180:270:0.6);
\draw [color=white, name path=hline-3] (-2,-0.2)--(-1.5,-0.2);
\draw [name intersections={of=circle-1 and hline-3}] (intersection-1)--(-2.8,-0.2);
\draw [dotted,thick](2.2,-0.4) arc(180:220:6.5 and 5.5);
\draw (2.2,-0.4) arc(180:205:6.5 and 5.5);
\draw (2.2,0.4)--(2.2,-0.4);
\draw (2.2,0.4) arc(0:180:2.55 and 2);
\draw (-2.9,0.4) arc(180:270:0.4);
\draw [color=white, name path=hline-2] (-2,0)--(-1.5,0);
\draw [dashed,name intersections={of=circle-1 and hline-2}] (intersection-1)--(-2.5,0);
\draw (3,-0.8) arc(180:270:6.5 and 4.5);
\draw (9.5,-5.3) arc(270:360:6.5 and 5.3);
\draw (16,0) arc(0:90:5.5 and 3.4);
\draw (10.5,3.4) arc(90:180:5.5 and 3);
\draw (5,0.4)--(5,-0.4);
\draw [dotted,thick](5,-0.4) arc(180:220:4.5 and 3);
\draw (5,-0.4) arc(180:200:4.5 and 3);
\draw (4,-0.8) arc(180:270:5.5 and 4);
\draw (9.5,-4.8) arc(270:360:5.5 and 4.8);
\draw (15,0) arc(0:90:4.5 and 2.9);
\draw (10.5,2.9) arc(90:180:4.5 and 2.5);
\draw (6,0.4)--(6,-0.4);
\draw [dotted,thick](6,-0.4) arc(180:220:3.5 and 2.5);
\draw (6,-0.4) arc(180:200:3.5 and 2.5);
\draw [name path=circle-2] (11.5,0) circle(0.5);
\draw [color=white, name path=hline-4] (12,-0.2)--(11.5,-0.2);
\draw [name intersections={of=circle-2 and hline-4}] (intersection-1)--(12.2,-0.2);
\draw (12.4,-0.4) arc(0:90:0.2);
\draw (10,-0.4) arc(180:360:1.2 and 1.2);
\draw (10,0.2)--(10,-0.4);
\draw [dotted,thick](10,0.2) arc(180:120:2 and 1);
\draw (10,0.2) arc(180:140:2 and 1);
\draw [color=white, name path=hline-6] (12,0.2)--(11.5,0.2);
\draw [name intersections={of=circle-2 and hline-6}] (intersection-1)--(12.8,0.2);
\draw (13.4,-0.4) arc(0:90:0.6);
\draw (9,-0.4) arc(180:360:2.2 and 2);
\draw (9,0.2)--(9,-0.4);
\draw [dotted,thick](9,0.2) arc(180:120:3 and 1.5);
\draw (9,0.2) arc(180:140:3 and 1.5);
\draw (4,-0.8)--(4,0.4);
\draw (4,0.4) arc(180:90:6.5 and 3.5);
\draw (10.5,3.9) arc(90:0:6.5 and 3.9);
\draw (17,0) arc(360:270:8 and 5.8);
\draw (9,-5.8) arc(270:180:8 and 5);
\draw (-1,0)--(1.9,0);
\draw (2.5,0)--(6.5,0);
\draw [dotted,thick] (6.6,0)--(8.4,0);
\draw (8.5,0)--(11,0);
\draw (1.9,0) arc(180:360:0.3);
\draw [dotted,thick] (13.5,0)--(14.8,0);
\draw [dashed](0,-6)--(0,-3);
\draw [>-<, dashed] (0,-3)--(0,3);
\draw [dashed](0,3)--(0,4);
\draw [<->, dashed] (-0.5,0)--(0.5,0);
\draw [dashed] (-1.5,0) arc(270:360:1.2 and 4);
\draw [->,dashed] (-0.3,4) arc(360:330:1.2 and 4);
\draw [dashed] (-1.5,0) arc(90:0:1.2 and 6);
\draw [->,dashed] (-0.3,-6) arc(0:30:1.2 and 6);
\draw [dashed] (-3.5,3)--(-1.5,0);
\draw [->,dashed] (-3.5,3)--(-2.5,1.5);
\draw [dashed] (-3.5,-3)--(-1.5,0);
\draw [->,dashed] (-3.5,-3)--(-2.5,-1.5);
\draw [dashed] (-3,0)--(-1.5,0);
\draw [->,dashed] (-3.5,0)--(-2.9,0);
\fill [white] (-1.5,0) circle(0.5);
\draw [thick] (-1.5,0) circle(0.5);
\draw [dashed](2.2,-6)--(2.2,-2);
\draw [<->, dashed] (2.2,-3)--(2.2,3);
\draw [dashed](2.2,2)--(2.2,4);
\draw [>-<, dashed] (1.7,0)--(3.7,0);
\draw [dashed] (0.3,4) arc(180:360:0.8 and 3);
\draw [->,dashed] (0.3,4) arc(180:210:0.8 and 3);
\draw [->,dashed] (0.3,4) arc(180:330:0.8 and 3);
\draw [dashed] (0.3,-6) arc(180:0:0.8 and 5);
\draw [->,dashed] (0.3,-6) arc(180:150:0.8 and 5);
\draw [->,dashed] (0.3,-6) arc(180:30:0.8 and 5);
\draw [dashed] (11.5,0.3) arc (270:180:8.5 and 3.7);
\draw [<->,dashed] (11.5,0.3) arc (270:220:8.5 and 3.7);
\draw [dashed] (11.5,-0.3) arc (90:180:8.5 and 5.7);
\draw [<->,dashed] (11.5,-0.3) arc (90:140:8.5 and 5.7);
\draw [dashed] (11.5,0)--(15.5,4);
\draw [<->,dashed] (11.5,0)--(13.5,2);
\draw [dashed] (11.5,0)--(15.5,-4);
\draw [<->,dashed] (11.5,0)--(13.5,-2);
\fill [white] (11.5,0) circle(0.5);
\draw [thick] (11.5,0) circle(0.5);
\draw (1.8,1.2) node[draw,fill] { \ \ \ \ \ \ \ \ \it q \ \  \ \ \ \ \ \ };
\draw (4.2,1.2) node[draw,fill] {\it r};
\draw (2.6,-1.6) node[draw,fill] {\ \ \ \ \ \  \ \ \it r+q \ \ \ \ \ \ \ \ \ \ };
\draw (9.5,-0.6) node[draw,fill] { \ \ \it q \ \ };
\draw (5.5,-0.6) node[draw,fill] { \ \ \it q \ \ };
\draw (2.2,-0.3) circle (1.5pt) node[anchor=north west] {\it $\Lambda_{1}$};
\draw [decorate, decoration=zigzag](0,0) arc(180:0:1.1 and 0.7);
\draw (0,0) arc(180:0:1.1 and 0.5);
\end{tikzpicture}
\caption{\it Framed Heegaard diagram when q is odd}
\label{fig:odd framed Heegaard}
\end{figure}

\begin{figure}\centering
\begin{tikzpicture}[scale=0.7]\scriptsize
\draw (-1.5,0) circle(0.5);
\draw (-0.4,-0.8)--(-0.4,0.4);
\draw (-0.4,0.4) arc(0:180:1 and 1);
\draw (-2.4,0.4) arc(180:270:0.4);
\draw (6.5,0) circle(0.5);
\draw (3.6,0.8)--(3.6,-0.4);
\draw (3.6,-0.4) arc(180:360:1.9 and 1);
\draw (7.4,-0.4) arc(0:90:0.4);
\draw (-1,0)--(6,0);
\fill[black, opacity=1] (-0.4,0) circle (1.5pt) node[anchor=south west] {\it $S^{m_j-2}(\Lambda_{j+1})$};
\fill[black, opacity=1] (3.6,0) circle (1.5pt) node[anchor=south west] {\it $S^{m_j}(\Lambda_{j})$};
\draw[->] (-0.4,-0.4)--(-0.4,-0.5);
\draw[->] (3.6,-0.3)--(3.6,-0.35);
\end{tikzpicture}
\caption{\it Power of S changing when q is odd}
\label{fig:pattern-1 when q odd}
\end{figure}

\begin{figure}\centering
\begin{tikzpicture}[scale=0.7]\scriptsize
\draw (-1.5,0) circle(0.5);
\draw (-0.4,0.8)--(-0.4,-0.4);
\draw (-0.4,-0.4) arc(180:360:4.1 and 1.7);
\draw (7.8,-0.4) arc(0:90:2.3 and 2.4);
\draw (5.5,2) arc(90:180:2.3 and 1.6);
\draw (6.5,0) circle(0.5);
\draw (3.2,0.4)--(3.2,-0.8);
\draw (-1,0)--(6,0);
\fill[black, opacity=1] (-0.4,0) circle (1.5pt) node[anchor=south west] {\it $S^{m_j}(\Lambda_{j})$};
\fill[black, opacity=1] (3.2,0) circle (1.5pt) node[anchor=south west] {\it $S^{m_j}(\Lambda_{j+1})$};
\draw[->] (3.2,-0.4)--(3.2,-0.5);
\draw[->] (-0.4,-0.3)--(-0.4,-0.35);
\end{tikzpicture}
\caption{\it Power of S changing when q is odd}
\label{fig:pattern-2 when q odd}
\end{figure}

\begin{figure}
\begin{tikzpicture}[scale=0.65,fill=white]\scriptsize
\begin{scope}
\draw (-1.5,0) circle(0.5);
\draw (0.5,-0.8)--(0.5,0.4);
\draw (0.5,0.4) arc(0:180:1.5 and 1);
\draw (-2.5,0.4) arc(180:270:0.5 and 0.4);
\draw (6.5,0) circle(0.5);
\draw (3.6,0.8)--(3.6,-0.4);
\draw (3.6,-0.4) arc(180:360:1.9 and 1);
\draw (7.4,-0.4) arc(0:90:0.4);
\draw (-1,0)--(6,0);
\draw [decorate, decoration=zigzag](-0.5,0) arc(180:0:1.1 and 0.7);
\draw (-0.5,0) arc(180:0:1.1 and 0.5);
\draw [dashed](-0.5,-1)--(-0.5,-2);
\draw [>-<, dashed] (-0.5,-1)--(-0.5,1);
\draw [dashed](-0.5,2)--(-0.5,1);
\draw [<->, dashed] (-0.8,0)--(-0.2,0);
\draw [dashed](1.7,-1)--(1.7,-2);
\draw [<->, dashed] (1.7,-1)--(1.7,1);
\draw [dashed](1.7,2)--(1.7,1);
\draw [>-<, dashed] (1.4,0)--(2,0);
\draw [color=white,line width=1mm] (1.7,-0.1)--(1.7,-0.5);
\fill[black, opacity=1] (0.5,0) circle (1.5pt) node[anchor=north west] {\it $S^{m_j}T(\Lambda_{j+1})$};
\fill[black, opacity=1] (3.6,0) circle (1.5pt) node[anchor=north west] {\it $S^{m_j}(\Lambda_{j})$};
\end{scope}
\begin{scope}[xshift = 12cm]
\draw (-1.5,0) circle(0.5);
\draw (1.5,-0.8)--(1.5,0.4);
\draw (1.5,0.4) arc(0:180:2 and 1);
\draw (-2.5,0.4) arc(180:270:0.5 and 0.4);
\draw (7.5,0) circle(0.5);
\draw (4.6,0.8)--(4.6,-0.4);
\draw (4.6,-0.4) arc(180:360:1.9 and 1);
\draw (8.4,-0.4) arc(0:90:0.4);
\draw (-1,0)--(7,0);
\draw [decorate, decoration=zigzag](-0.5,0) arc(180:0:0.8 and 0.7);
\draw (-0.5,0) arc(180:0:0.8 and 0.5);
\draw [dashed](-0.5,-1)--(-0.5,-2);
\draw [>-<, dashed] (-0.5,-1)--(-0.5,1);
\draw [dashed](-0.5,2)--(-0.5,1);
\draw [<->, dashed] (-0.8,0)--(-0.2,0);
\draw [dashed](1.1,-1)--(1.1,-2);
\draw [<->, dashed] (1.1,-1)--(1.1,1);
\draw [dashed](1.1,2)--(1.1,1);
\draw [>-<, dashed] (0.8,0)--(1.4,0);
\fill[black, opacity=1] (1.5,0) circle (1.5pt) node[anchor=north west] {\it $S^{m_j-2}(\Lambda_{j+1})$};
\fill[black, opacity=1] (4.6,0) circle (1.5pt) node[anchor=north west] {\it $S^{m_j}(\Lambda_{j})$};
\end{scope}
\end{tikzpicture}
\caption{\it Power of S changing when q is odd: case 1 and 2}
\label{fig:when q odd case 1 and 2}
\end{figure}

\begin{figure}
\begin{tikzpicture}[scale=0.65,fill=white]\scriptsize
\begin{scope}
\draw (-1.5,0) circle(0.5);
\draw (0.5,-0.8)--(0.5,0.4);
\draw (0.5,0.4) arc(0:180:1.5 and 1);
\draw (-2.5,0.4) arc(180:270:0.5 and 0.4);
\draw (6.5,0) circle(0.5);
\draw (3.6,0.8)--(3.6,-0.4);
\draw (3.6,-0.4) arc(180:360:1.9 and 1);
\draw (7.4,-0.4) arc(0:90:0.4);
\draw (-1,0)--(6,0);
\draw [decorate, decoration=zigzag](-0.5,0) arc(180:0:3 and 0.7);
\draw (-0.5,0) arc(180:0:3 and 0.5);
\draw [dashed](-0.5,-1)--(-0.5,-2);
\draw [>-<, dashed] (-0.5,-1)--(-0.5,1);
\draw [dashed](-0.5,2)--(-0.5,1);
\draw [<->, dashed] (-0.8,0)--(-0.2,0);
\draw [dashed](5.5,-1)--(5.5,-2);
\draw [<->, dashed] (5.5,-1)--(5.5,1);
\draw [dashed](5.5,2)--(5.5,1);
\draw [>-<, dashed] (5.2,0)--(5.8,0);
\draw [color=white,line width=1mm] (5.5,-0.1)--(5.5,-0.5);
\fill[black, opacity=1] (0.5,0) circle (1.5pt) node[anchor=north west] {\it $S^{m_j}T(\Lambda_{j+1})$};
\fill[black, opacity=1] (3.6,0) circle (1.5pt) node[anchor=north west] {\it $S^{m_j}(\Lambda_{j})$};
\end{scope}
\begin{scope}[xshift = 12cm]
\draw (-1.5,0) circle(0.5);
\draw (0.4,-0.4)--(0.4,0.8);
\draw (2.6,0.4) arc(0:180:2.6 and 1);
\draw (-2.6,0.4) arc(180:270:0.6 and 0.4);
\draw (6.5,0) circle(0.5);
\draw (2.6,0.4)--(2.6,-0.8);
\draw (0.4,-0.4) arc(180:360:3.5 and 1);
\draw (7.4,-0.4) arc(0:90:0.4);
\draw (-1,0)--(6,0);
\draw [decorate, decoration=zigzag](-0.5,0) arc(180:0:3 and 0.7);
\draw (-0.5,0) arc(180:0:3 and 0.5);
\draw [dashed](-0.5,-1)--(-0.5,-2);
\draw [>-<, dashed] (-0.5,-1)--(-0.5,1);
\draw [dashed](-0.5,2)--(-0.5,1);
\draw [<->, dashed] (-0.8,0)--(-0.2,0);
\draw [dashed](5.5,-1)--(5.5,-2);
\draw [<->, dashed] (5.5,-1)--(5.5,1);
\draw [dashed](5.5,2)--(5.5,1);
\draw [>-<, dashed] (5.2,0)--(5.8,0);
\fill[black, opacity=1] (0.4,0) circle (1.5pt) node[anchor=north west] {\it $S^{m_j}(\Lambda_{j})$};
\fill[black, opacity=1] (2.6,0) circle (1.5pt) node[anchor=north west] {\it $S^{m_j}T(\Lambda_{j+1})$};
\end{scope}
\end{tikzpicture}
\caption{\it Power of S changing when q is odd: case 3 and 4}
\label{fig:when q odd case 3 and 4}
\end{figure}

\begin{figure}
\begin{tikzpicture}[scale=0.65,fill=white]\scriptsize
\begin{scope}
\draw (-1.5,0) circle(0.5);
\draw (0.4,-0.4)--(0.4,0.8);
\draw (3,0.4) arc(0:180:2.8 and 1);
\draw (-2.6,0.4) arc(180:270:0.6 and 0.4);
\draw (6.5,0) circle(0.5);
\draw (3,0.4)--(3,-0.8);
\draw (0.4,-0.4) arc(180:360:3.5 and 1);
\draw (7.4,-0.4) arc(0:90:0.4);
\draw (-1,0)--(6,0);
\draw [decorate, decoration=zigzag](-0.5,0) arc(180:0:1.5 and 0.7);
\draw (-0.5,0) arc(180:0:1.5 and 0.5);
\draw [dashed](-0.5,-1)--(-0.5,-2);
\draw [>-<, dashed] (-0.5,-1)--(-0.5,1);
\draw [dashed](-0.5,2)--(-0.5,1);
\draw [<->, dashed] (-0.8,0)--(-0.2,0);
\draw [dashed](2.5,-1)--(2.5,-2);
\draw [<->, dashed] (2.5,-1)--(2.5,1);
\draw [dashed](2.5,2)--(2.5,1);
\draw [>-<, dashed] (2.2,0)--(2.8,0);
\fill[black, opacity=1] (0.4,0) circle (1.5pt) node[anchor=north west] {\it $S^{m_j}(\Lambda_{j})$};
\fill[black, opacity=1] (3,0) circle (1.5pt) node[anchor=north west] {\it $S^{m_j-2}(\Lambda_{j+1})$};
\end{scope}
\begin{scope}[xshift = 12cm]
\draw (-1.5,0) circle(0.5);
\draw (1.6,-0.4)--(1.6,0.8);
\draw (1.6,-0.4) arc(180:360:3.4 and 1.5);
\draw (7.5,0) circle(0.5);
\draw (8.4,-0.4) arc(0:90:2.2 and 1.9);
\draw (6.2,1.5) arc(90:180:2.2 and 1.1);
\draw (4,0.4)--(4,-0.8);
\draw (-1,0)--(7,0);
\draw [decorate, decoration=zigzag](-0.5,0) arc(180:0:0.8 and 0.7);
\draw (-0.5,0) arc(180:0:0.8 and 0.5);
\draw [dashed](-0.5,-1)--(-0.5,-2);
\draw [>-<, dashed] (-0.5,-1)--(-0.5,1);
\draw [dashed](-0.5,2)--(-0.5,1);
\draw [<->, dashed] (-0.8,0)--(-0.2,0);
\draw [dashed](1.1,-1)--(1.1,-2);
\draw [<->, dashed] (1.1,-1)--(1.1,1);
\draw [dashed](1.1,2)--(1.1,1);
\draw [>-<, dashed] (0.8,0)--(1.4,0);
\fill[black, opacity=1] (1.6,0) circle (1.5pt) node[anchor=north west] {\it $S^{m_j}(\Lambda_{j})$};
\fill[black, opacity=1] (4,0) circle (1.5pt) node[anchor=north west] {\it $S^{m_j}(\Lambda_{j+1})$};
\end{scope}
\end{tikzpicture}
\caption{\it Power of S changing when q is odd: case 5 and 6}
\label{fig:when q odd case 5 and 6}
\end{figure}

\begin{figure}\centering
\begin{tikzpicture}[scale=0.65,fill=white]\scriptsize
\draw (-1.5,0) circle(0.5);
\draw (0.4,-0.4)--(0.4,0.8);
\draw (6.5,0) circle(0.5);
\draw (3,0.4)--(3,-0.8);
\draw (0.4,-0.4) arc(180:360:3.5 and 1.5);
\draw (7.4,-0.4) arc(0:90:2.2 and 1.9);
\draw (5.2,1.5) arc(90:180:2.2 and 1.1);
\draw (-1,0)--(6,0);
\draw [decorate, decoration=zigzag](-0.5,0) arc(180:0:1.5 and 0.7);
\draw (-0.5,0) arc(180:0:1.5 and 0.5);
\draw [dashed](-0.5,-1)--(-0.5,-2);
\draw [>-<, dashed] (-0.5,-1)--(-0.5,1);
\draw [dashed](-0.5,2)--(-0.5,1);
\draw [<->, dashed] (-0.8,0)--(-0.2,0);
\draw [dashed](2.5,-1)--(2.5,-2);
\draw [<->, dashed] (2.5,-1)--(2.5,1);
\draw [dashed](2.5,2)--(2.5,1);
\draw [>-<, dashed] (2.2,0)--(2.8,0);
\fill[black, opacity=1] (0.4,0) circle (1.5pt) node[anchor=north west] {\it $S^{m_j}(\Lambda_{j})$};
\fill[black, opacity=1] (3,0) circle (1.5pt) node[anchor=north west] {\it $S^{m_j}(\Lambda_{j+1})$};
\end{tikzpicture}
\caption{\it Power of S changing when q is odd: case 7}
\label{fig:when q odd case 7}
\end{figure}
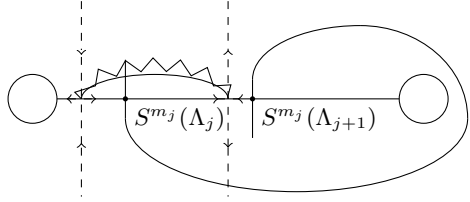

We choose a framed Heegaard diagram for $L(p,q)$ as shown in
Fig.~\ref{fig:odd framed Heegaard}. In this case, $k_1=1$. Let
us first analyze the pattern of powers of the antipode $S$ in the
product in Eq. ~\eqref{kuperberg invariant}. For $\Lambda_{k_1}$, which
is the starting point to do the multiplication along $c_u$, the
power of $S$ is
$2(\psi_l(k_1)-\psi_u(k_1))-\frac{1}{2}=2(-\frac{1}{4}-0)-\frac{1}{2}=-1$.

\begin{lemma}
The powers of $S$ change as shown in Fig.~\ref{fig:when q odd case
1 and 2} to Fig.~\ref{fig:when q odd case 7}. Namely, for $j=1,\ldots,p$,\\
(1) the power of $S$ from $\Lambda_{j}$ to the next factor $\Lambda_{j+1}$ decreases by $2$ when travelling along the 1-handle from right to left;\\
(2) the power of $S$ from $\Lambda_{j}$ to the next factor $\Lambda_{j+1}$ remains the same when travelling along the upper circle counterclockwise around the right disk.
\end{lemma}

\begin{proof}
Suppose the $j$-th term in the summation is
$S^{m_j}(\Lambda_{j})$, we calculate the difference
$m_{j+1}-m_j$ case by case, which is
$$m_{j+1}-m_j=2(\psi_l(\Lambda_{j+1})-\psi_l(\Lambda_{j}))-2(\psi_u(\Lambda_{j+1})-\psi_u(\Lambda_{j}))+2(\phi_u(\Lambda_{j+1})-\phi_u(\Lambda_{j}))$$
Here
$2(\phi_l(\Lambda_{j+1})-\phi_l(\Lambda_{j}))$ makes no contribution since the lower circle does not intersect with the twist fronts.
Note that $T=S^{-2}$ for factorizable Hopf algebras.\\
The patterns shown in Fig.~\ref{fig:pattern-1 when q odd} include five cases:\\
(1)\ \ This case is shown in the first picture in Fig.~\ref{fig:when
q odd case 1 and 2}.
$$m_{j+1}-m_j=2(0-(-\frac{1}{2}))-2(\frac{1}{2})+2(-1)=-2$$
(2)\ \ This case is shown in the second picture in
Fig.~\ref{fig:when q odd case 1 and 2}.
$$m_{j+1}-m_j=2((-\frac{1}{2})-(-\frac{1}{2}))-2(\frac{1}{2}+\frac{1}{2})+2(0)=-2$$
(3) and (4)\ \ These two cases are shown in the first picture in
Fig.~\ref{fig:when q odd case 3 and 4} and they share the same
calculation.
$$m_{j+1}-m_j=2(0-0)-2(-\frac{1}{2}+\frac{1}{2})+2(-1)=-2$$
(5)\ \ This case is shown in the first picture in Fig.~\ref{fig:when
q odd case 5 and 6}.
$$m_{j+1}-m_j=2(-\frac{1}{2}-0)-2(-\frac{1}{2}+\frac{1}{2}+\frac{1}{2})+2(0)=-2$$
The patterns shown in Fig.~\ref{fig:pattern-2 when q odd} include the following two cases:

(6)\ \ This case is shown in the second picture in
Fig.~\ref{fig:when q odd case 5 and 6}.
$$m_{j+1}-m_j=2((-\frac{1}{2})-(-\frac{1}{2}))-2(0)+2(0)=0$$
(7)\ \ This case is shown in the second picture in
Fig.~\ref{fig:when q odd case 7}.
$$m_{j+1}-m_j=2((-\frac{1}{2})-0)-2(-\frac{1}{2}))+2(0)=0$$
\end{proof}

Two more values to write down the Kuperberg invariant are
$\psi(c_l)$ and $\psi(c_u)$. It is easy to see that
$\psi(c_l)=-\frac{1}{2}$, and
$$\psi(c_u)=-\frac{1}{4}+(N_1-1)(\frac{1}{2}-\frac{1}{2})+(q-N_1)(\frac{1}{2}+\frac{1}{2})+\frac{1}{2}+\frac{1}{4}=\frac{q}{2}.$$
It follows that the power of $g$ is
$-\psi(c_u)+\frac{1}{2}=\frac{-q+1}{2}$.

Now we can write down the Kuperberg invariant
$Z_{Kup}(L(p,q),f,H)$
\begin{eqnarray*}
  Z_{Kup}&=&\prod\limits^q_{m=1}\prod\limits^{k_{m+1}-1}_{n=k_m}S^{-2m+1}(\Lambda_n)\\
  &=&\lambda\big(S^{-1}(\Lambda_{k_1})S^{-1}(\Lambda_{k_1+1})\cdots S^{-1}(\Lambda_{k_2-1})S^{-3}(\Lambda_{k_2})\cdots S^{-3}(\Lambda_{k_3-1})\\
  & &\cdots\cdots S^{-2q+1}(\Lambda_{k_q})\cdots S^{-2q+1}(\Lambda_{p})g^{\frac{-q+1}{2}}\big)\\
  &=&\lambda\big(S^{2q-1}(\Lambda_{k_1})S^{2q-1}(\Lambda_{k_1+1})\cdots S^{2q-1}(\Lambda_{k_2-1})S^{2q-3}(\Lambda_{k_2})\cdots S^{2q-3}(\Lambda_{k_3-1})\\
  & &\cdots\cdots S(\Lambda_{k_q})\cdots S(\Lambda_{p})g^{\frac{-q+1}{2}}\big)\\
  &=&\lambda\big(S^{2q-2}(\Lambda_{p})S^{2q-2}(\Lambda_{p-1})\cdots S^{2q-2}(\Lambda_{p+2-k_2})S^{2q-4}(\Lambda_{p+1-k_2})\cdots S^{2q-4}(\Lambda_{p+2-k_{3}})\\
  & &\cdots\cdots S^2(\Lambda_{p+1-k_{q-1}})\cdots S^2(\Lambda_{p+2-k_q})\Lambda_{p+1-k_q}\cdots \Lambda_{k_1+1}\Lambda_{k_1}g^{\frac{-q+1}{2}}\big)\\
  &=&\lambda\big(S^{2q-2}(\Lambda_{p})S^{2q-2}(\Lambda_{p-1})\cdots S^{2q-2}(\Lambda_{k_q})S^{2q-4}(\Lambda_{k_q-1})\cdots S^{2q-4}(\Lambda_{k_{q-1}})\\
  & &\cdots\cdots S^2(\Lambda_{k_3-1})\cdots S^2(\Lambda_{k_2})\Lambda_{k_2-1}\cdots \Lambda_{k_1+1}\Lambda_{k_1}g^{\frac{-q+1}{2}}\big)
\end{eqnarray*}
Here we have used the unimodular property that $S(\Lambda)=\Lambda$
and so $\sum S(\Lambda_{(p)})\otimes\cdots\otimes
S(\Lambda_{(p+1-j)})\otimes\cdots\otimes
S(\Lambda_{(1)})=\sum\Lambda_{(1)}\otimes\cdots\otimes
\Lambda_{(j)}\otimes\cdots\otimes\Lambda_{(p)}$ and by symmetry
$k_j+k_{q+2-j}=p+2$ for $j=2,\ldots,q$.

\subsubsection{$Z_{Kup}(L(p,q),f,H)$ when q is even}
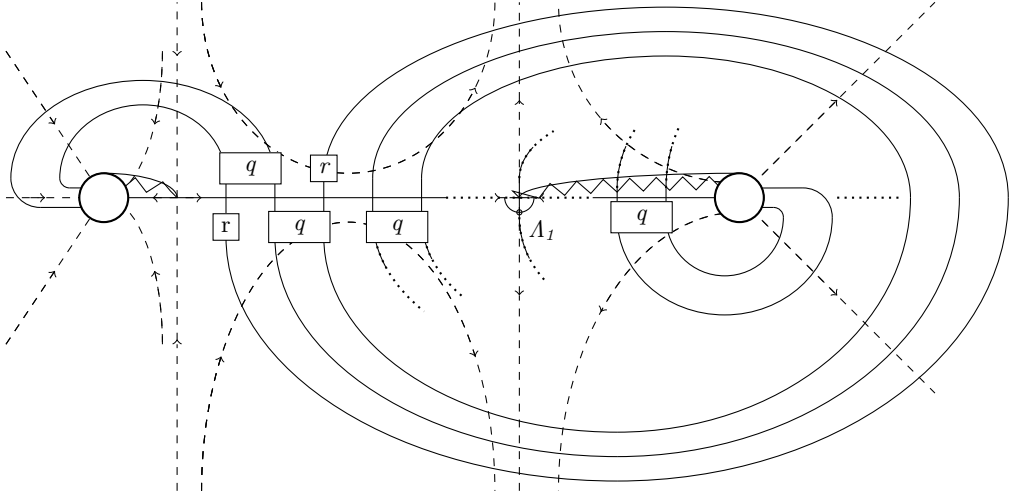
\begin{figure}\centering
\begin{tikzpicture}[scale=0.65,fill=white]\scriptsize
\draw [name path=circle-1] (-1.5,0) circle(0.5);
\draw (1,-0.8)--(1,0.4);
\draw (1,0.4) arc(0:180:1.7 and 1.5);
\draw (-2.4,0.4) arc(180:270:0.2);
\draw [color=white, name path=hline-1] (-2,0.2)--(-1.5,0.2);
\draw [name intersections={of=circle-1 and hline-1}] (intersection-1)--(-2.2,0.2);
\draw (2,-0.8)--(2,0.4);
\draw (2,0.4) arc(0:180:2.7 and 2);
\draw (-3.4,0.4) arc(180:270:0.6);
\draw [color=white, name path=hline-3] (-2,-0.2)--(-1.5,-0.2);
\draw [name intersections={of=circle-1 and hline-3}] (intersection-1)--(-2.8,-0.2);
\draw (2,-0.8) arc(180:270:7 and 4.5);
\draw (9,-5.3) arc(270:360:7 and 5.3);
\draw (16,0) arc(0:90:6 and 3.4);
\draw (10,3.4) arc(90:180:6 and 3);
\draw (4,0.4)--(4,-0.4);
\draw [dotted,thick](4,-0.4) arc(180:220:4.5 and 3);
\draw (4,-0.4) arc(180:200:4.5 and 3);
\draw (3,-0.8) arc(180:270:6 and 4);
\draw (9,-4.8) arc(270:360:6 and 4.8);
\draw (15,0) arc(0:90:5 and 2.9);
\draw (10,2.9) arc(90:180:5 and 2.5);
\draw (5,0.4)--(5,-0.4);
\draw [dotted,thick](5,-0.4) arc(180:220:3.5 and 2.5);
\draw (5,-0.4) arc(180:200:3.5 and 2.5);
\draw [name path=circle-2] (11.5,0) circle(0.5);
\draw [color=white, name path=hline-4] (12,-0.2)--(11.5,-0.2);
\draw [name intersections={of=circle-2 and hline-4}] (intersection-1)--(12.2,-0.2);
\draw (12.4,-0.4) arc(0:90:0.2);
\draw (10,-0.4) arc(180:360:1.2 and 1.2);
\draw (10,0.2)--(10,-0.4);
\draw [dotted,thick](10,0.2) arc(180:150:2 and 2.5);
\draw (10,0.2) arc(180:160:2 and 2.5);
\draw [color=white, name path=hline-6] (12,0.2)--(11.5,0.2);
\draw [name intersections={of=circle-2 and hline-6}] (intersection-1)--(12.8,0.2);
\draw (13.4,-0.4) arc(0:90:0.6);
\draw (9,-0.4) arc(180:360:2.2 and 2);
\draw (9,0.2)--(9,-0.4);
\draw [dotted,thick](9,0.2) arc(180:150:3 and 2.5);
\draw (9,0.2) arc(180:160:3 and 2.5);
\draw [dotted,thick](7,-0.4) arc(180:220:2.2 and 2);
\draw (7,-0.4) arc(180:205:2.2 and 2);
\draw (7,0.4)--(7,-0.4);
\draw [dotted,thick](7,0.4) arc(180:140:3 and 1.5);
\draw (7,0.4) arc(180:155:3 and 1.5);
\draw (3,-0.8)--(3,0.4);
\draw (3,0.4) arc(180:90:7 and 3.5);
\draw (10,3.9) arc(90:0:7 and 3.9);
\draw (17,0) arc(360:270:8 and 5.8);
\draw (9,-5.8) arc(270:180:8 and 5);
\draw (-1,0)--(5.5,0);
\draw [dotted,thick] (5.5,0)--(8.5,0);
\draw (8.5,0)--(11,0);
\draw (6.7,0) arc(180:360:0.3);
\draw [dotted,thick] (13.5,0)--(14.8,0);
\draw [decorate, decoration=zigzag](0,0) arc(0:90:1.5 and 0.3);
\draw (0,0) arc(0:90:1.5 and 0.5);
\draw [decorate, decoration=zigzag](7,0) arc(180:90:4.5 and 0.3);
\draw (7,0) arc(180:90:4.5 and 0.5);
\draw [dashed](0,-6)--(0,-3);
\draw [>-<, dashed] (0,-3)--(0,3);
\draw [dashed](0,3)--(0,4);
\draw [<->, dashed] (-0.5,0)--(0.5,0);
\draw [dashed] (-1.5,0) arc(270:360:1.2 and 3);
\draw [->,dashed] (-0.3,3) arc(360:330:1.2 and 3);
\draw [dashed] (-1.5,0) arc(90:0:1.2 and 3);
\draw [->,dashed] (-0.3,-3) arc(0:30:1.2 and 3);
\draw [dashed] (-3.5,3)--(-1.5,0);
\draw [->,dashed] (-3.5,3)--(-2.5,1.5);
\draw [dashed] (-3.5,-3)--(-1.5,0);
\draw [->,dashed] (-3.5,-3)--(-2.5,-1.5);
\draw [dashed] (-3,0)--(-1.5,0);
\draw [->,dashed] (-3.5,0)--(-2.7,0);
\fill [white] (-1.5,0) circle(0.5);
\draw [thick] (-1.5,0) circle(0.5);
\draw [dashed](7,-6)--(7,-2);
\draw [<->, dashed] (7,-2)--(7,2);
\draw [dashed](7,2)--(7,4);
\draw [>-<, dashed] (6.5,0)--(7.5,0);
\draw (6.8,0)--(7.2,0);
\draw [dashed] (0.5,4) arc(180:360:3 and 3.5);
\draw [->,dashed] (0.5,4) arc(180:210:3 and 3.5);
\draw [->,dashed] (0.5,4) arc(180:330:3 and 3.5);
\draw [dashed] (0.5,-6) arc(180:0:3 and 5.5);
\draw [->,dashed] (0.5,-6) arc(180:150:3 and 5.5);
\draw [->,dashed] (0.5,-6) arc(180:30:3 and 5.5);
\draw [dashed] (11.5,0.3) arc (270:180:3.7 and 3.7);
\draw [<->,dashed] (11.5,0.3) arc (270:220:3.7 and 3.7);
\draw [dashed] (11.5,-0.3) arc (90:180:3.7 and 5.7);
\draw [<->,dashed] (11.5,-0.3) arc (90:140:3.7 and 5.7);
\draw [dashed] (11.5,0)--(15.5,4);
\draw [<->,dashed] (11.5,0)--(13.5,2);
\draw [dashed] (11.5,0)--(15.5,-4);
\draw [<->,dashed] (11.5,0)--(13.5,-2);
\fill [white] (11.5,0) circle(0.5);
\draw [thick] (11.5,0) circle(0.5);
\draw (1.5,0.6) node[draw,fill] { \ \ \it q \ \ };
\draw (3,0.6) node[draw,fill] {\it r};
\draw (2.5,-0.6) node[draw,fill] { \ \ \it q \ \ };
\draw (1,-0.6) node[draw,fill] {r};
\draw (9.5,-0.4) node[draw,fill] { \ \ \it q \ \ };
\draw (4.5,-0.6) node[draw,fill] { \ \ \it q \ \ };
\draw (7,-0.3) circle (1.5pt) node[anchor=north west] {$\it \Lambda_{1}$};%
\end{tikzpicture}
\caption{\it Framed Heegaard diagram when q is even}
\label{fig:even framed Heegaard}
\end{figure}

\begin{figure}\centering
\begin{tikzpicture}[scale=0.7]\scriptsize
\draw (-1.5,0) circle(0.5);
\draw (-0.4,-0.8)--(-0.4,0.4);
\draw (-0.4,0.4) arc(0:180:1 and 1);
\draw (-2.4,0.4) arc(180:270:0.4);
\draw (6.5,0) circle(0.5);
\draw (2.6,0.8)--(2.6,-0.4);
\draw (2.6,-0.4) arc(180:360:2.4 and 1);
\draw (7.4,-0.4) arc(0:90:0.4);
\draw (-1,0)--(6,0);
\fill[black, opacity=1] (-0.4,0) circle (1.5pt) node[anchor=south west] {$S^{m_j}(\Lambda_{j+1})$};
\fill[black, opacity=1] (2.6,0) circle (1.5pt) node[anchor=south west] {$S^{m_j}(\Lambda_{j})$};
\end{tikzpicture}
\caption{\it Power of S changing when q is oven}
\label{fig:pattern-1 when q oven}
\end{figure}

\begin{figure}\centering
\begin{tikzpicture}[scale=0.7]\scriptsize
\draw (-1.5,0) circle(0.5);
\draw (-0.4,0.8)--(-0.4,-0.4);
\draw (-0.4,-0.4) arc(180:360:4.1 and 1.7);
\draw (7.8,-0.4) arc(0:90:2.7 and 2.4);
\draw (5.1,2) arc(90:180:2.7 and 1.6);
\draw (6.5,0) circle(0.5);
\draw (2.4,0.4)--(2.4,-0.4);
\draw (-1,0)--(6,0);
\fill[black, opacity=1] (-0.4,0) circle (1.5pt) node[anchor=south west] {$S^{m_j}(\Lambda_{j})$};
\fill[black, opacity=1] (2.4,0) circle (1.5pt) node[anchor=south west] {$S^{m_j+2}(\Lambda_{j+1})$};
\end{tikzpicture}
\caption{\it Power of S changing when q is oven}
\label{fig:pattern-2 when q oven}
\end{figure}

\begin{figure}
\begin{tikzpicture}[scale=0.65,fill=white]\scriptsize
\begin{scope}
\draw (-1.5,0) circle(0.5);
\draw (0.4,-0.8)--(0.4,0.4);
\draw (0.4,0.4) arc(0:180:1.5 and 1);
\draw (-2.6,0.4) arc(180:270:0.6 and 0.4);
\draw (6.5,0) circle(0.5);
\draw (3.6,0.8)--(3.6,-0.4);
\draw (3.6,-0.4) arc(180:360:1.9 and 1);
\draw (7.4,-0.4) arc(0:90:0.4);
\draw (-1,0)--(6,0);
\draw [dashed](0,-1)--(0,-2);
\draw [>-<, dashed] (0,-1)--(0,1);
\draw [dashed](0,2)--(0,1);
\draw [<->, dashed] (-0.3,0)--(0.3,0);
\draw [dashed](3,-1)--(3,-2);
\draw [<->, dashed] (3,-1)--(3,1);
\draw [dashed](3,2)--(3,1);
\draw [>-<, dashed] (2.7,0)--(3.2,0);
\draw [decorate, decoration=zigzag](0,0) arc(0:90:1.5 and 0.3);
\draw (0,0) arc(0:90:1.5 and 0.5);
\draw [decorate, decoration=zigzag](3,0) arc(180:90:3.5 and 0.3);
\draw (3,0) arc(180:90:3.5 and 0.5);
\fill [white] (-1.5,0) circle(0.5);
\draw [thick] (-1.5,0) circle(0.5);
\fill [white] (6.5,0) circle(0.5);
\draw [thick] (6.5,0) circle(0.5);
\fill[black, opacity=1] (0.4,0) circle (1.5pt) node[anchor=north west] {$S^{m_j}(\Lambda_{j+1})$};
\fill[black, opacity=1] (3.6,0) circle (1.5pt) node[anchor=north west] {$S^{m_j}(\Lambda_{j})$};
\end{scope}
\begin{scope}[xshift = 12cm]
\draw (-1.5,0) circle(0.5);
\draw (0.4,-0.8)--(0.4,0.4);
\draw (0.4,0.4) arc(0:180:1.5 and 1);
\draw (-2.6,0.4) arc(180:270:0.6 and 0.4);
\draw (6.5,0) circle(0.5);
\draw (3,0.8)--(3,-0.4);
\draw (3,-0.4) arc(180:360:2.2 and 1);
\draw (7.4,-0.4) arc(0:90:0.4);
\draw (-1,0)--(6,0);
\draw [dashed](0,-1)--(0,-2);
\draw [>-<, dashed] (0,-1)--(0,1);
\draw [dashed](0,2)--(0,1);
\draw [<->, dashed] (-0.3,0)--(0.3,0);
\draw [dashed](3.6,-1)--(3.6,-2);
\draw [<->, dashed] (3.6,-1)--(3.6,1);
\draw [dashed](3.6,2)--(3.6,1);
\draw [>-<, dashed] (3.3,0)--(3.9,0);
\draw [color=white,line width=1mm] (3.6,-0.1)--(3.6,-0.5);
\draw [decorate, decoration=zigzag](0,0) arc(0:90:1.5 and 0.3);
\draw (0,0) arc(0:90:1.5 and 0.5);
\draw [decorate, decoration=zigzag](3.6,0) arc(180:90:2.9 and 0.3);
\draw (3.6,0) arc(180:90:2.9 and 0.5);
\fill [white] (-1.5,0) circle(0.5);
\draw [thick] (-1.5,0) circle(0.5);
\fill [white] (6.5,0) circle(0.5);
\draw [thick] (6.5,0) circle(0.5);
\fill[black, opacity=1] (0.4,0) circle (1.5pt) node[anchor=north west] {$S^{m_j}(\Lambda_{j+1})$};
\fill[black, opacity=1] (3,0) circle (1.5pt) node[anchor=north west] {$S^{m_j}(\Lambda_{j})$};
\end{scope}
\end{tikzpicture}
\caption{\it Power of S changing when q is even: case 1 and 2}
\label{fig:when q even case 1 and 2}
\end{figure}

\begin{figure}
\begin{tikzpicture}[scale=0.65,fill=white]\scriptsize
\begin{scope}
\draw (-1.5,0) circle(0.5);
\draw (0.4,-0.4)--(0.4,0.8);
\draw (7,0) circle(0.5);
\draw (3,0.4)--(3,-0.8);
\draw (0.4,-0.4) arc(180:360:3.8 and 1.5);
\draw (8,-0.4) arc(0:90:2.5 and 1.9);
\draw (5.5,1.5) arc(90:180:2.5 and 1.1);
\draw (-1,0)--(6.5,0);
\draw [dashed](0,-1)--(0,-2);
\draw [>-<, dashed] (0,-1)--(0,1);
\draw [dashed](0,2)--(0,1);
\draw [<->, dashed] (-0.3,0)--(0.3,0);
\draw [dashed](2.5,-1)--(2.5,-2);
\draw [<->, dashed] (2.5,-1)--(2.5,1);
\draw [dashed](2.5,2)--(2.5,1);
\draw [>-<, dashed] (2.2,0)--(2.8,0);
\draw [decorate, decoration=zigzag](0,0) arc(0:90:1.5 and 0.3);
\draw (0,0) arc(0:90:1.5 and 0.5);
\draw [decorate, decoration=zigzag](2.5,0) arc(180:90:4.5 and 0.3);
\draw (2.5,0) arc(180:90:4.5 and 0.5);
\fill [white] (-1.5,0) circle(0.5);
\draw [thick] (-1.5,0) circle(0.5);
\fill [white] (7,0) circle(0.5);
\draw [thick] (7,0) circle(0.5);
\fill[black, opacity=1] (0.4,0) circle (1.5pt) node[anchor=north west] {$S^{m_j}(\Lambda_{j})$};
\fill[black, opacity=1] (3,0) circle (1.5pt) node[anchor=north west] {$S^{m_j}T^{-1}(\Lambda_{j+1})$};
\end{scope}
\begin{scope}[xshift = 12cm]
\draw (-1.5,0) circle(0.5);
\draw (0.4,-0.4)--(0.4,0.8);
\draw (7,0) circle(0.5);
\draw (3,0.4)--(3,-0.8);
\draw (0.4,-0.4) arc(180:360:3.8 and 1.5);
\draw (8,-0.4) arc(0:90:2.5 and 1.9);
\draw (5.5,1.5) arc(90:180:2.5 and 1.1);
\draw (-1,0)--(6.5,0);
\draw [dashed](0,-1)--(0,-2);
\draw [>-<, dashed] (0,-1)--(0,1);
\draw [dashed](0,2)--(0,1);
\draw [<->, dashed] (-0.3,0)--(0.3,0);
\draw [dashed](3.5,-1)--(3.5,-2);
\draw [<->, dashed] (3.5,-1)--(3.5,1);
\draw [dashed](3.5,2)--(3.5,1);
\draw [>-<, dashed] (3.2,0)--(3.8,0);
\draw [color=white,line width=1mm] (3.5,-0.1)--(3.5,-0.5);
\draw [decorate, decoration=zigzag](0,0) arc(0:90:1.5 and 0.3);
\draw (0,0) arc(0:90:1.5 and 0.5);
\draw [decorate, decoration=zigzag](3.5,0) arc(180:90:3.5 and 0.3);
\draw (3.5,0) arc(180:90:3.5 and 0.5);
\fill [white] (-1.5,0) circle(0.5);
\draw [thick] (-1.5,0) circle(0.5);
\fill [white] (7,0) circle(0.5);
\draw [thick] (7,0) circle(0.5);
\fill[black, opacity=1] (0.4,0) circle (1.5pt) node[anchor=north west] {$S^{m_j}(\Lambda_{j})$};
\fill[black, opacity=1] (3,0) circle (1.5pt) node[anchor=north west] {$S^{m_j+2}(\Lambda_{j+1})$};
\end{scope}
\end{tikzpicture}
\caption{\it Power of S changing when q is even: case 3 and 4}
\label{fig:when q even case 3 and 4}
\end{figure}

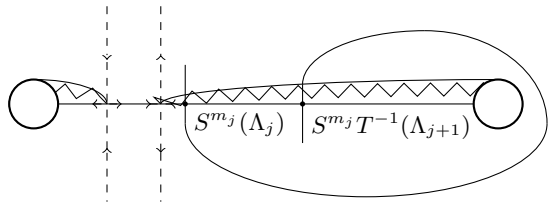
\begin{figure}\centering
\begin{tikzpicture}[scale=0.65,fill=white]\scriptsize
\draw (-1.5,0) circle(0.5);
\draw (1.6,-0.4)--(1.6,0.8);
\draw (1.6,-0.4) arc(180:360:3.7 and 1.5);
\draw (8,0) circle(0.5);
\draw (9,-0.4) arc(0:90:2.5 and 1.9);
\draw (6.5,1.5) arc(90:180:2.5 and 1.1);
\draw (4,0.4)--(4,-0.8);
\draw (-1,0)--(7.5,0);
\draw [dashed](0,-1)--(0,-2);
\draw [>-<, dashed] (0,-1)--(0,1);
\draw [dashed](0,2)--(0,1);
\draw [<->, dashed] (-0.3,0)--(0.3,0);
\draw [dashed](1.1,-1)--(1.1,-2);
\draw [<->, dashed] (1.1,-1)--(1.1,1);
\draw [dashed](1.1,2)--(1.1,1);
\draw [>-<, dashed] (0.8,0)--(1.4,0);
\draw [decorate, decoration=zigzag](0,0) arc(0:90:1.5 and 0.3);
\draw (0,0) arc(0:90:1.5 and 0.5);
\draw [decorate, decoration=zigzag](1.1,0) arc(180:90:6.9 and 0.3);
\draw (1.1,0) arc(180:90:6.9 and 0.5);
\fill [white] (-1.5,0) circle(0.5);
\draw [thick] (-1.5,0) circle(0.5);
\fill [white] (8,0) circle(0.5);
\draw [thick] (8,0) circle(0.5);
\fill[black, opacity=1] (1.6,0) circle (1.5pt) node[anchor=north west] {$S^{m_j}(\Lambda_{j})$};
\fill[black, opacity=1] (4,0) circle (1.5pt) node[anchor=north west] {$S^{m_j}T^{-1}(\Lambda_{j+1})$};
\end{tikzpicture}
\caption{\it Power of S changing when q is even: case 5}
\label{fig:when q even case 5}
\end{figure}

Fig.~\ref{fig:even framed Heegaard} is a framed Heegaard diagram for
$L(p,q)$ when $q$ is even. In this case, the twist front is
different from the case when $q$ is odd, so the pattern of the
power changes of $S$ is also different.
\begin{lemma}
As shown in Fig.~\ref{fig:pattern-1 when q oven} and Fig.~\ref{fig:pattern-2
when q oven}, for $j=1,\ldots,p$,\\
(1) the power of $S$ from $\Lambda_{j}$ to the next factor $\Lambda_{j+1}$ remains the same when travelling along the upper circle counterclockwise around the right disk;\\
(2) the power of $S$ from $\Lambda_{j}$ to the next factor $\Lambda_{j+1}$ increases by $2$ when travelling along the 1-handle from right to left.
\end{lemma}

\begin{proof}
Similar to the $q$ odd case, we calculate the change $m_{j+1}-m_j$ of the powers
of $S$ case by case shown from Fig.~\ref{fig:when q even case 1 and 2} to
Fig.~\ref{fig:when q even case 5}.\\
(1) as shown in the first picture in Fig.~\ref{fig:when q even case
1 and 2},
$$m_{j+1}-m_j=2(0-(-\frac{1}{2}))-2(\frac{1}{2})+2(0)=0.$$
(2) as shown in the second picture in Fig.~\ref{fig:when q even case
1 and 2},
$$m_{j+1}-m_j=2(0-0)-2(-\frac{1}{2}+\frac{1}{2})+2(0)=0.$$
(3) as shown in the first picture in Fig.~\ref{fig:when q even case
3 and 4},
$$m_{j+1}-m_j=2((-\frac{1}{2})-0)-2(-\frac{1}{2})+2(1)=2.$$
(4) as shown in the second picture in Fig.~\ref{fig:when q even case
3 and 4},
$$m_{j+1}-m_j=2(0-0)-2(-\frac{1}{2}-\frac{1}{2})+2(0)=2.$$
(5) as shown in the second picture in Fig.~\ref{fig:when q even case
5},
$$m_{j+1}-m_j=2((-\frac{1}{2})-(-\frac{1}{2}))-2(0-0)+2(1)=2$$
\end{proof}

For $\psi(c_u)$, we have
$$\psi(c_u)=-\frac{1}{4}+(N_1-1-q)(-\frac{1}{2}-\frac{1}{2})-\frac{1}{4}=\frac{-p+q}{2}$$
and then the power of $g$ is $-\psi(c_u)+\frac{1}{2}=\frac{p-q+1}{2}$.

Thus the Kuperberg invariant is
\begin{eqnarray*}
  Z_{Kup}&=&\prod\limits^q_{m=0}\prod\limits^{k_{m+1}-1}_{n=k_m}S^{2n-2m-3}(\Lambda_n)\\
  &=&\lambda\big(S^{-1}(\Lambda_{1})S(\Lambda_{2})\cdots S^{2k_1-5}(\Lambda_{k_1-1})S^{2k_1-5}(\Lambda_{k_1})S^{2k_1-3}(\Lambda_{k_1+1})\cdots S^{2k_2-7}(\Lambda_{k_2-1})\\
  & &\cdots\cdots S^{2k_q-2q-3}(\Lambda_{k_q})S^{2k_q-2q-1}(\Lambda_{k_q+1})\cdots S^{2p-2q-3}(\Lambda_{p})g^{\frac{p-q+1}{2}}\big)\\
  &=&\lambda\big(\Lambda_{p}S^2(\Lambda_{p-1})\cdots S^{2k_1-4}(\Lambda_{p+2-k_1})S^{2k_1-4}(\Lambda_{p+1-k_1})S^{2k_1-2}(\Lambda_{p-k_1})\cdots S^{2k_2-6}(\Lambda_{p-k_2+1})\\
  & &\cdots\cdots S^{2k_q-2q-2}(\Lambda_{p+1-k_q})S^{2k_q-2q}(\Lambda_{p-k_q})\cdots S^{2p-2q-2}(\Lambda_{1})g^{\frac{p-q+1}{2}}\big)\\
  &=&\lambda\big(S^{-2p+2q+2}(\Lambda_{p})S^{-2p+2q+4}(\Lambda_{p-1})\cdots S^{-2p+2q+2k_1-2}(\Lambda_{p+2-k_1})\\
  & &S^{-2p+2q+2k_1-2}(\Lambda_{p+1-k_1})S^{-2p+2q+2k_1}(\Lambda_{p-k_1})\cdots S^{-2p+2q+2k_2-4}(\Lambda_{p-k_2+1})\\
  & &\cdots\cdots S^{-2p+2k_q}(\Lambda_{p+2-k_q})S^{-2p+2k_q}(\Lambda_{p+1-k_q})S^{-2p+2k_q+2}(\Lambda_{p-k_q})\cdots S^{-2}(\Lambda_{2})\Lambda_{1}g^{\frac{p-q+1}{2}}\big)\\
  &=&\lambda\big(S^{-2p+2q+2}(\Lambda_{p})S^{-2p+2q+4}(\Lambda_{p-1})\cdots S^{-2k_q+2q+2}(\Lambda_{k_q})\\
  & &S^{-2k_q+2q+2}(\Lambda_{k_q-1})S^{-2k_q+2q+4}(\Lambda_{k_q-2})\cdots S^{-2k_{q-1}+2q}(\Lambda_{k_{q-1}-1})\\
  & &\cdots\cdots S^{-2k_1+4}(\Lambda_{k_1})S^{-2k_1+4}(\Lambda_{k_1-1})S^{-2k_1+6}(\Lambda_{k_1-2})\cdots S^{-2}(\Lambda_{2})\Lambda_{1}g^{\frac{p-q+1}{2}}\big)\\
\end{eqnarray*}
Here we have used that $\sum
S^{-1}(\Lambda_{(p)})\otimes\cdots\otimes
S^{-1}(\Lambda_{(p+1-j)})\otimes\cdots\otimes
S^{-1}(\Lambda_{(1)})=\sum\Lambda_{(1)}\otimes\cdots\otimes
\Lambda_{(j)}\otimes\cdots\otimes\Lambda_{(p)}$ and $k_j+k_{q+1-j}=p+2$ for $j=1,\ldots,q$.

\subsection{$Z_{Henn}(L(p,q)\#\overline{L(p,q)},H)$}
We use the chain-mail link $L$ in Fig.~\ref{fig:Chain-mail L(p,q)}
to calculate the Hennings invariant for $L(p,q)\#\overline{L(p,q)}$.
Since the signature $\sigma(L)$ of the framing matrix of the chain-mail link is zero, so the normalization factor
$$[\lambda(\theta)\lambda(\theta^{-1})]^{-\frac{c(L)}{2}}[\lambda(\theta)/\lambda(\theta^{-1})]^{-\frac{\sigma(L)}{2}}=1$$
because $\lambda(\theta)\lambda(\theta^{-1})=1$ for
a factorizable ribbon Hopf algebra (see~\cite{CW2}). Hence it is sufficient to find the
link invariant $TR(L,H)$. First, the contribution of the lower circle $c_l$ to the Hennings invariant
is equivalent to decorating the upper circle $c_u$ with cointegrals. That is
\begin{lemma}
\[
\begin{tikzpicture}[scale=0.8]
\draw [name path=ellipse] (0,0) ellipse (3cm and 1cm);
\draw (1,-1.5)--(1,1.5);
\draw (2,-1.5)--(2,1.5);
\draw (-1,-1.5)--(-1,1.5);
\draw (-2,-1.5)--(-2,1.5);
\draw (3,0) [color=white,line width=2mm] arc (0:180:3cm and 1cm);
\draw (3,0) arc (0:180:3cm and 1cm);
\draw [color=white,line width=2mm](1,-1.5)--(1,0);
\draw [color=white,line width=2mm](2,-1.5)--(2,0);
\draw [color=white,line width=2mm](-1,-1.5)--(-1,0);
\draw [color=white,line width=2mm](-2,-1.5)--(-2,0);
\draw (1,-1.5)--(1,0);
\draw (2,-1.5)--(2,0);
\draw (-1,-1.5)--(-1,0);
\draw (-2,-1.5)--(-2,0);
\draw [color=white,name path=vline](0.5,-1.5)--(0.5,1.5);
\draw [color=white, line width=2mm, name intersections={of=ellipse and vline}] (intersection-1) arc (80:100:3cm and 1cm);
\draw [dashed, name intersections={of=ellipse and vline}] (intersection-1) arc (80:100:3cm and 1cm);
\draw [color=white, line width=2mm, name intersections={of=ellipse and vline}] (intersection-2) arc (280:260:3cm and 1cm);
\draw [dashed, name intersections={of=ellipse and vline}] (intersection-2) arc (280:260:3cm and 1cm);
\draw (-3,0) circle (0.1pt) node[anchor=east] {$c_l$};
\node (=) at (4.2,0.2) {$=$};
\begin{scope}[xshift = 9cm]
\draw (-100pt,30pt)--(-100pt,-30pt);
\draw (-50pt,30pt)--(-50pt,-30pt);
\draw (100pt,30pt)--(100pt,-30pt);
\draw (50pt,30pt)--(50pt,-30pt);
\draw[dotted] (-30pt,10pt)--(30pt,10pt);
\fill (-100pt,-10pt) circle (1.5pt) node[anchor=west] {$\Lambda_{(1)}$};
\fill (-50pt,-10pt) circle (1.5pt) node[anchor=west] {$\Lambda_{(2)}$};
\fill (50pt,-10pt) circle (1.5pt) node[anchor=west] {$\Lambda_{(p-1)}$};
\fill (100pt,-10pt) circle (1.5pt) node[anchor=west] {$\Lambda_{(p)}$};
\end{scope}
\end{tikzpicture}
\]
\end{lemma}

\begin{proof}
By Kauffman and Radford's algorithm, we label the crossings with components of the ${R}$-matrix $R=\sum s\otimes t$, where $\sum s^1\otimes
t^1=\cdots=\sum s^{2p}\otimes t^{2p}$ are copies of the $R$-matrix, and obtain the immersed diagram in Fig.~\ref{fig:immersion lower circle}.
Then we can separate the circle from the rest and push all the decorated elements
to one side as shown in Fig.~\ref{fig:seperated lower circle}. This
diagram gives us the following tensor element. Here the last
equality results from lemma 2.
\begin{eqnarray*}
& &\sum \lambda(t^{2p}t^{2p-1}\cdots t^{p+2}t^{p+1}s^ps^{p-1}\cdots s^2s^1) s^{2p}t^1\otimes s^{2p-1}t^2\otimes\cdots\otimes s^{p+2}t^{p-1}\otimes s^{p+1}t^p\\
&=&\sum \lambda_{(1)}(t^{2p})\lambda_{(2p)}(s^1) s^{2p}t^1\otimes\lambda_{(2)}(t^{2p-1})\lambda_{(2p-2)}(s^2) s^{2p-1}t^2\otimes\cdots\\
& &\otimes\lambda_{(p-1)}(t^{p+2})\lambda_{(p+1)}(s^{p-1}) s^{p+2}t^{p-1}\otimes\lambda_{(p)}(t^{p+1}s^p) s^{p+1}t^p\\
&=&f_{R^{\tau}}(\lambda_{(1)})f_R(\lambda_{(2p)})\otimes f_{R^{\tau}}(\lambda_{(2)})f_R(\lambda_{(2p-2)})\otimes\ldots\otimes f_{R^{\tau}}(\lambda_{(p-1)})f_R(\lambda_{(p+1)})\otimes f_{R^{\tau}R}(\lambda_{(p)})\\
&=&\sum\Lambda_{(1)}\otimes\ldots\otimes\Lambda_{(p)}
\end{eqnarray*}

\begin{figure}\centering
\begin{tikzpicture}[scale=1,rotate=-45]\scriptsize
\draw[color=white,name path=ellipse1] (0,0) ellipse (130pt and 50pt);
\draw[name path=ellipse2] (0,0) ellipse (120pt and 40pt);
\draw[name path=vline11] (-80pt,60pt)--(-80pt,-60pt);
\draw[name path=vline12] (-40pt,60pt)--(-40pt,-60pt);
\draw[name path=vline13] (80pt,60pt)--(80pt,-60pt);
\draw[name path=vline14] (40pt,60pt)--(40pt,-60pt);
\draw[color=white,name path=vline21] (-90pt,60pt)--(-90pt,-60pt);
\draw[color=white,name path=vline22] (-50pt,60pt)--(-50pt,-60pt);
\draw[color=white,name path=vline23] (70pt,60pt)--(70pt,-60pt);
\draw[color=white,name path=vline24] (30pt,60pt)--(30pt,-60pt);
\draw[color=white,name path=vline31] (-70pt,60pt)--(-70pt,-60pt);
\draw[color=white,name path=vline32] (-30pt,60pt)--(-30pt,-60pt);
\draw[color=white,name path=vline33] (90pt,60pt)--(90pt,-60pt);
\draw[color=white,name path=vline34] (50pt,60pt)--(50pt,-60pt);
\fill[black, opacity=1, name intersections={of= ellipse2 and vline21}] (intersection-1) circle (1.5pt) node[anchor=south] {$s^1$};
\fill[black, opacity=1, name intersections={of= ellipse1 and vline11}] (intersection-1) circle (1.5pt) node[anchor=south] {$t^1$};
\fill[black, opacity=1, name intersections={of= ellipse2 and vline22}] (intersection-1) circle (1.5pt) node[anchor=south] {$s^2$};
\fill[black, opacity=1, name intersections={of= ellipse1 and vline12}] (intersection-1) circle (1.5pt) node[anchor=south] {$t^2$};
\fill[black, opacity=1, name intersections={of= ellipse2 and vline23}] (intersection-1) circle (1.5pt) node[anchor=east] {$s^p$};
\fill[black, opacity=1, name intersections={of= ellipse1 and vline13}] (intersection-1) circle (1.5pt) node[anchor=west] {$t^p$};
\fill[black, opacity=1, name intersections={of= ellipse2 and vline24}] (intersection-1) circle (1.5pt) node[anchor=east] {$s^{p-1}$};
\fill[black, opacity=1, name intersections={of= ellipse1 and vline14}] (intersection-1) circle (1.5pt) node[anchor=west] {$t^{p-1}$};
\fill[black, opacity=1, name intersections={of= ellipse1 and vline11}] (intersection-2) circle (1.5pt) node[anchor=east] {$s^{2p}$};
\fill[black, opacity=1, name intersections={of= ellipse2 and vline31}] (intersection-2) circle (1.5pt) node[anchor=west] {$S^{-1}(t^{2p})$};
\fill[black, opacity=1, name intersections={of= ellipse1 and vline12}] (intersection-2) circle (1.5pt) node[anchor=east] {$s^{2p-1}$};
\fill[black, opacity=1, name intersections={of= ellipse2 and vline32}] (intersection-2) circle (1.5pt) node[anchor=west] {$S^{-1}(t^{2p-1})$};
\fill[black, opacity=1, name intersections={of= ellipse1 and vline13}] (intersection-2) circle (1.5pt) node[anchor=east] {$s^{p+1}$};
\fill[black, opacity=1, name intersections={of= ellipse2 and vline33}] (intersection-2) circle (1.5pt) node[anchor=south west] {$S^{-1}(t^{p+1})$};
\fill[black, opacity=1, name intersections={of= ellipse1 and vline14}] (intersection-2) circle (1.5pt) node[anchor=east] {$s^{p+2}$};
\fill[black, opacity=1, name intersections={of= ellipse2 and vline34}] (intersection-2) circle (1.5pt) node[anchor=south west] {$S^{-1}(t^{p+2})$};
\draw[color=white,line width=2mm] (20pt,39.5pt) arc (80:100:120pt and 40pt);
\draw[dashed] (20pt,39.5pt) arc (80:100:120pt and 40pt);
\draw[color=white,line width=2mm] (20pt,-39.5pt) arc (-80:-100:120pt and 40pt);
\draw[dashed] (20pt,-39.5pt) arc (-80:-100:120pt and 40pt);
\end{tikzpicture}
\caption{\it Immersion diagram of lower circle in chain-mail link} \label{fig:immersion lower circle}
\end{figure}

\begin{figure}\centering
\begin{tikzpicture}[scale=0.8]\scriptsize
\draw[name path=ellipse] (200pt,0pt) ellipse (40pt and 120pt);
\draw (-100pt,30pt)--(-100pt,-30pt);
\draw (-50pt,30pt)--(-50pt,-30pt);
\draw (100pt,30pt)--(100pt,-30pt);
\draw (50pt,30pt)--(50pt,-30pt);
\draw[dotted] (-30pt,10pt)--(30pt,10pt);
\fill (-100pt,-10pt) circle (1.5pt) node[anchor=east] {$s^{2p}t^1$};
\fill (-50pt,-10pt) circle (1.5pt) node[anchor=east] {$s^{2p-1}t^2$};
\fill (50pt,-10pt) circle (1.5pt) node[anchor=east] {$s^{p+2}t^{p-1}$};
\fill (100pt,-10pt) circle (1.5pt) node[anchor=east] {$s^{p+1}t^p$};
\draw[color=white,name path=hline1] (200pt,110pt)--(250pt,110pt);
\draw[color=white,name path=hline2] (200pt,90pt)--(250pt,90pt);
\draw[color=white,name path=hline31,line width=1mm] (200pt,75pt)--(250pt,75pt);
\draw[color=white,name path=hline32,line width=1mm] (200pt,70pt)--(250pt,70pt);
\draw[color=white,name path=hline33,line width=1mm] (200pt,65pt)--(250pt,65pt);
\draw[color=white,name path=hline34,line width=1mm] (200pt,60pt)--(250pt,60pt);
\draw[color=white,name path=hline35,line width=1mm] (200pt,55pt)--(250pt,55pt);
\draw[color=white,name path=hline4] (200pt,50pt)--(250pt,50pt);
\draw[color=white,name path=hline5] (200pt,30pt)--(250pt,30pt);
\draw[color=white,name path=hline6] (200pt,10pt)--(250pt,10pt);
\draw[color=white,name path=hline7] (200pt,-10pt)--(250pt,-10pt);
\draw[color=white,name path=hline8] (200pt,-30pt)--(250pt,-30pt);
\draw[color=white,name path=hline9] (200pt,-50pt)--(250pt,-50pt);
\draw[color=white,name path=hline101,line width=1mm] (200pt,-75pt)--(250pt,-75pt);
\draw[color=white,name path=hline102,line width=1mm] (200pt,-70pt)--(250pt,-70pt);
\draw[color=white,name path=hline103,line width=1mm] (200pt,-65pt)--(250pt,-65pt);
\draw[color=white,name path=hline104,line width=1mm] (200pt,-60pt)--(250pt,-60pt);
\draw[color=white,name path=hline105,line width=1mm] (200pt,-55pt)--(250pt,-55pt);
\draw[color=white,name path=hline11] (200pt,-90pt)--(250pt,-90pt);
\draw[color=white,name path=hline12] (200pt,-110pt)--(250pt,-110pt);
\fill[black, opacity=1, name intersections={of= ellipse and hline1}] (intersection-1) circle (1.5pt) node[anchor=west] {$s^1$};
\fill[black, opacity=1, name intersections={of= ellipse and hline2}] (intersection-1) circle (1.5pt) node[anchor=west] {$s^2$};
\fill[black, opacity=1, name intersections={of= ellipse and hline5}] (intersection-1) circle (1.5pt) node[anchor=west] {$s^{p-1}$};
\fill[black, opacity=1, name intersections={of= ellipse and hline6}] (intersection-1) circle (1.5pt) node[anchor=west] {$s^p$};
\fill[black, opacity=1, name intersections={of= ellipse and hline7}] (intersection-1) circle (1.5pt) node[anchor=west] {$t^{p+1}$};
\fill[black, opacity=1, name intersections={of= ellipse and hline8}] (intersection-1) circle (1.5pt) node[anchor=west] {$t^{p+2}$};
\fill[black, opacity=1, name intersections={of= ellipse and hline11}] (intersection-1) circle (1.5pt) node[anchor=west] {$t^{2p-1}$};
\fill[black, opacity=1, name intersections={of= ellipse and hline12}] (intersection-1) circle (1.5pt) node[anchor=west] {$t^{2p}$};
\end{tikzpicture}
\caption{\it Immersion diagram of lower circle in chain-mail link} \label{fig:seperated lower circle}
\end{figure}
\end{proof}

\begin{figure}\centering
\begin{tikzpicture}
\draw (-80pt,30pt)--(-80pt,-30pt);
\draw (-40pt,30pt)--(-40pt,-30pt);
\draw (80pt,30pt)--(80pt,-30pt);
\draw (40pt,30pt)--(40pt,-30pt);
\draw[color=white,line width=2mm] (-100pt,10pt)--(100pt,10pt);
\draw (-100pt,10pt)--(-25pt,10pt);
\draw[dotted] (-25pt,10pt)--(25pt,10pt);
\draw (25pt,10pt)--(100pt,10pt);
\fill (-80pt,-20pt) circle (1.5pt) node[anchor=west] {$\Lambda_{(1)}$};
\fill (-40pt,-20pt) circle (1.5pt) node[anchor=west] {$\Lambda_{(2)}$};
\fill (40pt,-20pt) circle (1.5pt) node[anchor=west] {$\Lambda_{(p-1)}$};
\fill (80pt,-20pt) circle (1.5pt) node[anchor=west] {$\Lambda_{(p)}$};
\end{tikzpicture}
\caption{\it Typical crossing in upper circle}
\label{fig:typical crossing}
\end{figure}
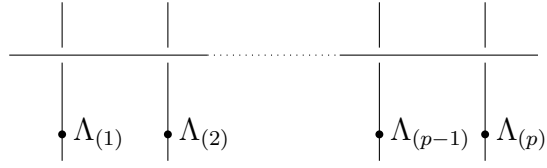

The next step is to resolve the crossings where the upper circle
$c_u$ crosses itself. A typical crossing in the chain-mail is shown
in Fig.~\ref{fig:typical crossing}.

\begin{lemma}
  \[
\begin{tikzpicture}[scale=0.9]\footnotesize
\draw (-80pt,30pt)--(-80pt,-30pt);
\draw (-40pt,30pt)--(-40pt,-30pt);
\draw (80pt,30pt)--(80pt,-30pt);
\draw (40pt,30pt)--(40pt,-30pt);
\draw[color=white,line width=2mm] (-100pt,10pt)--(100pt,10pt);
\draw (-100pt,10pt)--(-25pt,10pt);
\draw[dotted] (-25pt,10pt)--(25pt,10pt);
\draw (25pt,10pt)--(100pt,10pt);
\fill (-80pt,-20pt) circle (1.5pt) node[anchor=west] {$\Lambda_{(1)}$};
\fill (-40pt,-20pt) circle (1.5pt) node[anchor=west] {$\Lambda_{(2)}$};
\fill (40pt,-20pt) circle (1.5pt) node[anchor=west] {$\Lambda_{(p-1)}$};
\fill (80pt,-20pt) circle (1.5pt) node[anchor=west] {$\Lambda_{(p)}$};
\node (=) at (4,0.2) {$=$};
\begin{scope}[xshift = 8cm]
\draw (-80pt,30pt)--(-80pt,-30pt);
\draw (-40pt,30pt)--(-40pt,-30pt);
\draw (80pt,30pt)--(80pt,-30pt);
\draw (40pt,30pt)--(40pt,-30pt);
\draw (-100pt,10pt)--(-25pt,10pt);
\draw[dotted] (-25pt,10pt)--(25pt,10pt);
\draw (25pt,10pt)--(100pt,10pt);
\fill (-80pt,-20pt) circle (1.5pt) node[anchor=west] {$\Lambda_{(1)}$};
\fill (-40pt,-20pt) circle (1.5pt) node[anchor=west] {$\Lambda_{(2)}$};
\fill (40pt,-20pt) circle (1.5pt) node[anchor=west] {$\Lambda_{(p-1)}$};
\fill (80pt,-20pt) circle (1.5pt) node[anchor=west] {$\Lambda_{(p)}$};
\end{scope}
\end{tikzpicture}
  \]
\end{lemma}

\begin{proof}
The corresponding immersed diagram is\\
\begin{tikzpicture}[scale=0.8]
\begin{scope}[rotate=-45]
\draw (-80pt,30pt)--(-80pt,-60pt);
\draw (-40pt,30pt)--(-40pt,-50pt);
\draw (80pt,30pt)--(80pt,-30pt);
\draw (40pt,30pt)--(40pt,-40pt);
\draw (-100pt,10pt)--(-25pt,10pt);
\draw[dotted] (-25pt,10pt)--(25pt,10pt);
\draw (25pt,10pt)--(100pt,10pt);
\fill (-80pt,-20pt) circle (1.5pt) node[anchor=east] {$\Lambda_{(1)}$};
\fill (-40pt,-20pt) circle (1.5pt) node[anchor=east] {$\Lambda_{(2)}$};
\fill (40pt,-20pt) circle (1.5pt) node[anchor=east] {$\Lambda_{(p-1)}$};
\fill (80pt,-20pt) circle (1.5pt) node[anchor=east] {$\Lambda_{(p)}$};
\fill (-90pt,10pt) circle (1.5pt) node[anchor=east] {$s^1$};
\fill (-80pt,20pt) circle (1.5pt) node[anchor=west] {$t^1$};
\fill (-50pt,10pt) circle (1.5pt) node[anchor=east] {$s^2$};
\fill (-40pt,20pt) circle (1.5pt) node[anchor=west] {$t^2$};
\fill (30pt,10pt) circle (1.5pt) node[anchor=east] {$s^{p-1}$};
\fill (40pt,20pt) circle (1.5pt) node[anchor=west] {$t^{p-1}$};
\fill (70pt,10pt) circle (1.5pt) node[anchor=east] {$s^p$};
\fill (80pt,20pt) circle (1.5pt) node[anchor=west] {$t^p$};
\end{scope}
\node (=) at (4,0.2) {$=$};
\begin{scope}[xshift = 10cm,rotate=-45]
\draw (-80pt,30pt)--(-80pt,-60pt);
\draw (-40pt,30pt)--(-40pt,-50pt);
\draw (80pt,30pt)--(80pt,-30pt);
\draw (40pt,30pt)--(40pt,-40pt);
\draw (-100pt,10pt)--(-25pt,10pt);
\draw[dotted] (-25pt,10pt)--(25pt,10pt);
\draw (25pt,10pt)--(100pt,10pt);
\fill (-80pt,-20pt) circle (1.5pt) node[anchor=east] {$\Lambda_{(1)}$};
\fill (-40pt,-20pt) circle (1.5pt) node[anchor=east] {$\Lambda_{(2)}$};
\fill (40pt,-20pt) circle (1.5pt) node[anchor=east] {$\Lambda_{(p-1)}$};
\fill (80pt,-20pt) circle (1.5pt) node[anchor=east] {$\Lambda_{(p)}$};
\fill (-80pt,20pt) circle (1.5pt) node[anchor=west] {$t^1$};
\fill (-40pt,20pt) circle (1.5pt) node[anchor=west] {$t^2$};
\fill (40pt,20pt) circle (1.5pt) node[anchor=west] {$t^{p-1}$};
\fill (80pt,20pt) circle (1.5pt) node[anchor=west] {$t^p$};
\fill (-90pt,10pt) circle (1.5pt) node[anchor=east] {$s^ps^{p-1}\cdots s^2s^1$};
\end{scope}
\end{tikzpicture}

From this immersed diagram, we obtain the following tensor element
\begin{eqnarray*}
& &\sum (s^ps^{p-1}\cdots s^2s^1)\otimes\Lambda_{(1)}t^1\otimes\Lambda_{(2)}t^2\otimes\cdots\otimes\Lambda_{(p-1)}t^{p-1}\otimes\Lambda_{(p)}t^p\\
&=&\sum s\otimes\Lambda_{(1)}t^{(1)}\otimes\Lambda_{(2)}t^{(2)}\otimes\cdots\otimes\Lambda_{(p-1)}t^{(p-1)}\otimes\Lambda_{(p)}t^{(p)}\\
&=&\sum \varepsilon(t)s\otimes\Lambda_{(1)}\otimes\Lambda_{(2)}\otimes\cdots\otimes\Lambda_{(p-1)}\otimes\Lambda_{(p)}\\
&=&\sum 1\otimes\Lambda_{(1)}\otimes\Lambda_{(2)}\otimes\cdots\otimes\Lambda_{(p-1)}\otimes\Lambda_{(p)}\\
&=&\sum \Lambda_{(1)}\otimes\Lambda_{(2)}\otimes\cdots\otimes\Lambda_{(p-1)}\otimes\Lambda_{(p)}
\end{eqnarray*}
Here we have used the property of the $R$-matrix that
$(\varepsilon\otimes id)(R)=1$ and
$(id\otimes\Delta^{(p-1)})(R)=R_{1,p+1}\cdots R_{12}$ which comes
from $(id\otimes\Delta)(R)=R_{13}R_{12}$. \\
\end{proof}

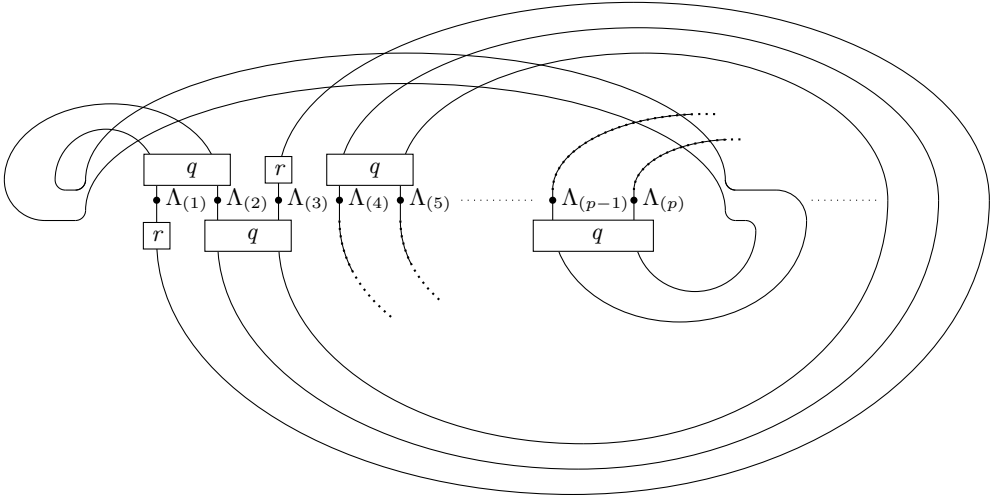
\begin{figure}
\begin{tikzpicture}[scale=0.675,fill=white]\scriptsize
\draw (-0.4,-0.8)--(-0.4,0.4);
\draw (-0.4,0.4) arc(0:180:1 and 1);
\draw (-2.4,0.4) arc(180:270:0.2);
\draw (-2,0.2)--(-2.2,0.2);
\draw (0.8,-0.8)--(0.8,0.4);
\draw (0.8,0.4) arc(0:180:2.1 and 1.5);
\draw (-3.4,0.4) arc(180:270:0.8);
\draw (-2,-0.4)--(-2.6,-0.4);
\draw (0.8,-0.8) arc(180:270:7.1 and 4.5);
\draw (7.9,-5.3) arc(270:360:7.1 and 5.3);
\draw (15,0) arc(0:90:5.9 and 3.4);
\draw (9.1,3.4) arc(90:180:5.9 and 3);
\draw (3.2,0.4)--(3.2,-0.4);
\draw [dotted,thick](3.2,-0.4) arc(180:220:4.5 and 3);
\draw (3.2,-0.4) arc(180:200:4.5 and 3);
\draw (2,-0.8) arc(180:270:6 and 4);
\draw (8,-4.8) arc(270:360:6 and 4.8);
\draw (14,0) arc(0:90:4.8 and 2.9);
\draw (9.2,2.9) arc(90:180:4.8 and 2.5);
\draw (4.4,0.4)--(4.4,-0.4);
\draw [dotted,thick](4.4,-0.4) arc(180:220:3.5 and 2.5);
\draw (4.4,-0.4) arc(180:200:3.5 and 2.5);
\draw (11,-0.4)--(11.2,-0.4);
\draw (11.4,-0.6) arc(0:90:0.2);
\draw (9,-0.6) arc(180:360:1.2 and 1.2);
\draw (9,0.2)--(9,-0.6);
\draw [dotted,thick](9,0.2) arc(180:85:2 and 1);
\draw (9,0.2) arc(180:95:2 and 1);
\draw (11,0.2)--(11.8,0.2);
\draw (12.4,-0.4) arc(0:90:0.6);
\draw (7.4,-0.4) arc(180:360:2.5 and 2);
\draw (7.4,0.2)--(7.4,-0.4);
\draw [dotted,thick](7.4,0.2) arc(180:85:3 and 1.5);
\draw (7.4,0.2) arc(180:95:3 and 1.5);
\draw (2,-0.8)--(2,0.4);
\draw (2,0.4) arc(180:90:7 and 3.5);
\draw (9,3.9) arc(90:0:7 and 3.9);
\draw (16,0) arc(360:270:8.2 and 5.8);
\draw (7.8,-5.8) arc(270:180:8.2 and 5);
\draw [dotted] (5.6,0)--(7,0);
\draw [dotted] (12.5,0)--(13.8,0);
\draw (0.2,0.6) node[draw,fill] { \ \ \ \ \it q\ \ \ \ };
\draw (2,0.6) node[draw,fill] {\it r~~};
\draw (1.4,-0.7) node[draw,fill] { \ \ \ \ \it q\ \ \ \ };
\draw (-0.4,-0.7) node[draw,fill] {\it r~~};
\draw (8.2,-0.7) node[draw,fill] { \ \ \ \ \ \ \it q\ \ \ \  \ \  };
\draw (3.8,0.6) node[draw,fill] { \ \ \ \ \it q\ \ \ \ };
\draw (-2,0.2) arc(270:360:0.2);
\draw (10.8,0.4) arc(0:180:6.3 and 2.5);
\draw (11,0.2) arc(270:180:0.2);
\draw (-2,-0.4) arc(270:360:0.2);
\draw (10.8,-0.2) arc(0:180:6.3 and 2.5);
\draw (11,0.2) arc(270:180:0.2);
\draw (11,-0.4) arc(270:180:0.2);
\draw (-2,0.2) arc(270:360:0.2);
\fill[black, opacity=1] (-0.4,0) circle (2pt) node[anchor=west] {$\Lambda_{(1)}$};
\fill[black, opacity=1] (0.8,0) circle (2pt) node[anchor=west] {$\Lambda_{(2)}$};
\fill[black, opacity=1] (2,0) circle (2pt) node[anchor=west] {$\Lambda_{(3)}$};
\fill[black, opacity=1] (3.2,0) circle (2pt) node[anchor=west] {$\Lambda_{(4)}$};
\fill[black, opacity=1] (4.4,0) circle (2pt) node[anchor=west] {$\Lambda_{(5)}$};
\fill[black, opacity=1] (7.4,0) circle (2pt) node[anchor=west] {$\Lambda_{(p-1)}$};
\fill[black, opacity=1] (9,0) circle (2pt) node[anchor=west] {$\Lambda_{(p)}$};
\end{tikzpicture}
\caption{\it Chain-mail link decorated with cointegrals}
\label{fig:chain-mail with cointegrals}
\end{figure}

The last step is to push all the labels $\Lambda_{(i)}$'s in
Fig.~\ref{fig:chain-mail with cointegrals} to where
$\Lambda_{(N_1)}$ is located and then do the evaluation by $\lambda$
to get the Hennings invariants. 

In the example of $L(2,1)$, we push all the labels to $\Lambda_{(1)}$, then 
$$Z_{Henn}(L(2,1)\#\overline{L(2,1)})=\lambda\big(S^{-2}(\Lambda_{(2)})\Lambda_{(1)}G^{2}\big).$$
To compare with $Z_{Kup}(L(2,1))$, we use $G^{-1}S^2(x)=xG^{-1}$ for $x\in H$ and $\Lambda G^{-1}=\Lambda$, then obtain
\begin{eqnarray*}
  Z_{Henn}(L(2,1)\#\overline{L(2,1)})&=&\lambda\big(S^{-2}(\Lambda_{(2)})\Lambda_{(1)}G^{2}\big)\\
  &=&\lambda\big(S^{-2}(\Lambda_{(2)})G^{-1}\Lambda_{(1)}G^{-1}G^{2}\big)\\
  &=&\lambda\big(G^{-1}\Lambda_{(2)}\Lambda_{(1)}G^{}\big)\\
  &=&\lambda\big(\Lambda_{(2)}\Lambda_{(1)}\big)
\end{eqnarray*}
That is equal to $Z_{Kup}(L(2,1))$.

For $L(5,2)$, we push all the labels to $\Lambda_{(4)}$, then 
$$Z_{Henn}(L(5,2)\#\overline{L(5,2)})=\lambda\big(S^{-4}(\Lambda_{(2)})S^{-4}(\Lambda_{(5)})S^{-2}(\Lambda_{(3)})\Lambda_{(1)}\Lambda_{(4)}G^{4}\big).$$
which is the same as $Z_{Kup}(L(5,2))$.

The general case will be done 
according to whether $q$ is odd or even.

\subsubsection{$Z_{Henn}(L(p,q)\#\overline{L(p,q)},H)$ when $q$ is odd}
In this case, we push all the labels $\Lambda_{(i)}$'s to
$\Lambda_{(N_0)}$ along the upper circle and write down the
following equality for $Z_{Henn}(L(p,q)\#\overline{L(p,q)},H)$.

\begin{eqnarray*}
  Z_{Henn}&=&\lambda\big(S^{-2p+2q}(\Lambda_{p})\cdots S^{-2k_{q}+2q}(\Lambda_{k_q})\cdots\cdots S^{-2k_3+6}(\Lambda_{k_3-1})\cdots S^{-2k_2+4}(\Lambda_{k_2})\\
  & &  S^{-2k_2+4}(\Lambda_{k_2-1})\cdots S^{-2}(\Lambda_{k_1+1})\Lambda_{k_1}G^{p-q+1}\big)\\
  &=&\lambda\big(S^{-2p+2q}(\Lambda_{p})G^{-1}\cdots S^{-2k_{q}+2q}(\Lambda_{k_q})G^{-1}\cdots\cdots S^{-2k_3+6}(\Lambda_{k_3-1})G^{-1}\cdots\\
  & &S^{-2k_2+4}(\Lambda_{k_2})G^{-1}S^{-2k_2+4}(\Lambda_{k_2-1})G^{-1}\cdots S^{-2}(\Lambda_{k_1+1})G^{-1}\Lambda_{k_1}G^{-1}G^{p-q+1}\big)\\
  &=&\lambda\big(S^{2q-2}(\Lambda_{p})\cdots S^{2q-2}(\Lambda_{k_q})\cdots\cdots S^{2}(\Lambda_{k_3-1})\cdots S^{2}(\Lambda_{k_2})\\
  & &\Lambda_{k_2-1}\cdots \Lambda_{k_1+1}\Lambda_{k_1}G^{-q+1}\big)
\end{eqnarray*}
Here we have used the fact that $G$ is grouplike and
$G^{-1}S^2(x)=xG^{-1}$ and $\Lambda G^{-1}=\Lambda$. Note that $G^2=g$. Hence we obtain
that when $q$ is odd, $Z_{Kup}(L(p,q),f,H)=Z_{Henn}(L(p,q)\#\overline{L(p,q)},H)$.

\subsubsection{$Z_{Henn}(L(p,q)\#\overline{L(p,q)},H)$ when $q$ is even}
Now, we push all the labels $\Lambda_{(i)}$'s to
$\Lambda_{1}$ along the upper circle and obtain
\begin{eqnarray*}
Z_{Henn}&=&\lambda\big(S^{-2p+2q+2}(\Lambda_{p})S^{-2p+2q+4}(\Lambda_{p-1})\cdots S^{-2k_q+2q+2}(\Lambda_{k_q})\\
  & &S^{-2k_q+2q+2}(\Lambda_{k_q-1})S^{-2k_q+2q+4}(\Lambda_{k_q-2})\cdots S^{-2k_{q-1}+2q}(\Lambda_{k_{q-1}-1})\\
  & &\cdots\cdots S^{-2k_1+4}(\Lambda_{k_1})S^{-2k_1+4}(\Lambda_{k_1-1})\cdots S^{-2}(\Lambda_{2})\Lambda_{1}G^{p-q+1}\big)
\end{eqnarray*}
Thus $Z_{Kup}(L(p,q),f,H)=Z_{Henn}(L(p,q)\#\overline{L(p,q)},H)$ when $q$ is even.

\subsection{Remarks on the general case}

It is natural to conjecture that the relation $|Z_{Kup}(M,f,H)|=|Z_{Henn}(M,H)|^2$ always holds for any closed oriented $3$-manifold $M$ and factorizable finite dimensional ribbon Hopf algebra $H$.
For a general closed oriented $3$-manifold $M$, inspired by the result in \cite{CKS}, one strategy to prove the conjecture is to divide the problem into two cases:

\begin{enumerate}

\item If $H_1(M,\Q)\neq 0$, then both invariants of $M$ are $0$;

\item If $H_1(M,\Q)=0$, then a similar comparison can be carried out.

\end{enumerate}

Unfortunately, the choice of a suitable framing for the Kuperberg invariant which allows a direct comparison with the Hennings invariant is extremely hard to come by.  Even for the lens spaces, we are lucky to find the suitable framings.  Some other choices that we tried led to expressions that were hard to compare the two invariants.  Ideas on fermionic TQFTs in \cite {GWW} might be relevant for a conceptual approach to the conjecture.


\begin{thebibliography}{AA}

\bibitem[BK]{BK}B. Balsam and A. Kirillov Jr., Turaev-Viro invariants as an extended TQFT, arXiv:1004.1533.

\bibitem[BW]{BW} J. Barrett and B. Westbury,
\textit{The equality of 3-manifold invariants}, Math. Proc. Cambridge Philos. Soc.\ \textbf{118}, (1995), 503-510.

\bibitem[CKS]{CKS} Q. Chen, S. Kuppum, and P. Srinivasan,
\textit{On the relation between the WRT invariant and the Hennings invariant}, Math. Proc. Cambridge Philos. Soc.\ \textbf{146}, (2009), 151-163.

\bibitem[CW1]{CW1} M. Cohen and S. Westreich,
\textit{Some interrelations between Hopf algebras and their duals}, J. Algebra \ \textbf{283}, 1 (2005), 42-62.

\bibitem[CW2]{CW2} M. Cohen and S. Westreich,
\textit{Fourier transforms for Hopf algebras}, Contemp. Math. \ \textbf{433}, (2007), 115-133.

\bibitem[GWW]{GWW}Z.-G. Gu, Z. Wang, X.-G. Wen,
 A classification of 2D fermionic and bosonic topological orders, arXiv:1010.1517.

\bibitem[H]{H} M. A. Hennings,
\textit{Invariants of links and 3-manifolds obtained form Hopf algebras}, J. London Math. Soc.\ \textbf{54}, 2 (1996), 594-624.

\bibitem[Ker]{Kerler1}T.\ Kerler, Genealogy of non-perturbative quantum-invariants of $3$-manifolds: the surgical family. Geometry and physics (Aarhus, 1995), 503?547, Lecture Notes in Pure and Appl. Math., 184, Dekker, New York, 1997


\bibitem[KL]{nonsemisimple}T.\ Kerler; V. V. Lyubashenko, Non-Semisimple Topological Quantum Field Theories for 3-Manifolds with Corners,
Lecture Notes in Mathematics, 1765. Springer-Verlag, Berlin, 2001. vi+379 pp.

\bibitem[KR]{KR} L. H. Kauffman and D. E. Radford,
\textit{Invariants of 3-manifolds derived form finite dimensional Hopf algebras}, J. Knot Theory Ramifications\ \textbf{4}, 1 (1995), 131-162.

\bibitem[KR2]{KR2} L. H. Kauffman and D. E. Radford,
\textit{A Necessary and Sufficient Condition for a Finite-Dimensional Drinfel'd Double to Be a Ribbon Hopf Algebra}, J. Algebra\ \textbf{159}, 1 (1993), 98-114.

\bibitem[Ku]{Ku} G. Kuperberg,
\textit{Noninvolutory Hopf algebras and 3-manifold invariants}, Duke Math. J.\ \textbf{84}, 1 (1996), 83-129.


\bibitem[Ra1]{Ra1} D. Radford,
\textit{The trace function and Hopf algebras}, J. Algebra\ \textbf{163}, 3 (1994), 583-622.

\bibitem[Ra2]{Ra2} D. Radford,
\textit{On Kauffman's knot invariants arising from finite-dimensional Hopf algebras}, In Advances in Hopf algebras (Chicago, IL, 1992), pages 205--266.  Dekker, New York, 1994.

\bibitem[Ro]{Ro} J. Roberts,
\textit{Skein theory and Turaev-Viro invariants}, Topology\ \textbf{34}, 4 (1995), 771-787.


\bibitem[TV1]{TV1} V. Turaev and O. Viro,
\emph{State sum invariants of $3$-manifolds and quantum $6j$-symbols.}
Topology \textbf{31} (1992), no. 4, 865--902.

\bibitem[TV2]{TV2} V. Turaev, A. Virelizier, On two approaches to 3-dimensional TQFTs, arXiv:1006.3501


\end{thebibliography}
\end{document}